\documentclass[11pt, a4paper]{article}
\usepackage[utf8]{inputenc} 
\usepackage[a4paper, margin=2.5cm]{geometry}
\usepackage{parskip}
\usepackage{comment}
\usepackage{enumitem}
\usepackage{multirow}
\usepackage{cite}

\usepackage{mathtools} 
\usepackage{amsthm}
\usepackage{amssymb}
\usepackage{bbm}
\usepackage{mathrsfs}
\usepackage{accents}
\usepackage{stackengine}
\usepackage{scalerel}

\mathtoolsset{showonlyrefs}

\usepackage{xcolor}
\usepackage[export]{adjustbox}
\usepackage{tikz,pgfplots}

\usetikzlibrary{
  arrows,
  arrows.meta,
  bending,
  calc,
  decorations.markings,
  decorations.text,
  intersections,
  patterns,
  spy
}

\pgfplotsset{compat=1.18}

\usepackage{algorithm}
\usepackage{algorithmic}

\definecolor{hanblue}{rgb}{0.27, 0.42, 0.81}
\definecolor{mordantred19}{rgb}{0.68, 0.05, 0.0}
\definecolor{myred}{rgb}{0.8, 0.0, 0.0}
\definecolor{mygreen}{rgb}{0.0, 0.5, 0.0}

\usepackage[
  colorlinks,
  citecolor=hanblue,
  linkcolor=mordantred19
]{hyperref}

\DeclareMathAlphabet{\mathpgoth}{OT1}{pgoth}{m}{n}

\DeclareFontFamily{U}{mathx}{}
\DeclareFontShape{U}{mathx}{m}{n}{<-> mathx10}{}
\DeclareSymbolFont{mathx}{U}{mathx}{m}{n}
\DeclareMathAccent{\widehat}{0}{mathx}{"70}
\DeclareMathAccent{\widecheck}{0}{mathx}{"71}

\tikzset{
  mid arrow/.style={
    postaction={
      decorate,
      decoration={
        markings,
        mark=at position 0.5 with {\arrow{Stealth[length=4pt,width=5pt]}}
      }
    }
  }
}

\newcommand{\calA}{\mathcal{A}}

\newcommand{\calG}{\mathcal{G}}

\newcommand{\calP}{\mathcal{P}}

\newcommand{\calR}{\mathcal{R}}

\newcommand{\calS}{\mathcal{S}}

\newcommand{\calL}{\mathcal{L}}

\newcommand{\bC}{\mathbb{C}}
\newcommand{\bH}{\mathbb{H}}          
\newcommand{\bG}{\mathbb{G}}

\newcommand{\R}{\mathbb{R}}
\newcommand{\Z}{\mathbb{Z}}

\newcommand{\rH}{\mathring{\mathbb{H}}} 
\newcommand{\C}{\mathbb{C}}

\renewcommand{\iota}{\mathrm{i}}

\newcommand{\Hb}{\mathbb{H}}

\renewcommand{\epsilon}{\varepsilon}

\newcommand{\dd}{\,\mathrm{d}}

\def\mint{\mkern12mu\hbox{\vrule height4pt depth-3.2pt width5pt}\mkern-16.5mu\int}

\newcommand\restr[2]{{
  \left.\kern-\nulldelimiterspace
  #1
  \vphantom{\big|}
  \right|_{#2}
}}

\theoremstyle{plain}
\newtheorem{theorem}{Theorem}[section]

\newtheorem{lemma}[theorem]{Lemma}
\newtheorem{proposition}[theorem]{Proposition}
\newtheorem{corollary}[theorem]{Corollary}

\newtheorem{assumption}{Assumption}
\newtheorem{remark}[theorem]{Remark}

\theoremstyle{definition}
\newtheorem{definition}[theorem]{Definition}

\theoremstyle{remark}

\title{Phase-equivariant signal analysis via optimal transport in the Heisenberg group}

\author{Giacomo Cristinelli$^\ast$, Gautam Pai$^\ast$, Christoph Brune$^\ast$, Jos\'e A. Iglesias\thanks{Department of Applied Mathematics, University of Twente, 7500AE Enschede, The Netherlands \newline (\texttt{g.cristinelli@utwente.nl, g.pai{@}utwente.nl, c.brune{@}utwente.nl, jose.iglesias{@}utwente.nl})}}
\date{}

\begin{document}

\maketitle

\begin{abstract}
    We introduce a phase-sensitive time-frequency representation of complex signals as probability measures on the reduced Heisenberg group. The construction lifts the short-time Fourier transform by retaining both its magnitude and phase, and is equivariant under the Schr\"odinger action on the signal. Endowing the Heisenberg group with a left-invariant (subRiemannian) metric therefore yields Wasserstein comparisons that are invariant under common delay, frequency and phase transformations. We establish reconstruction and inverse stability properties of the lift, with errors controlled by the discarded low-amplitude coefficients and by Wasserstein perturbations. We illustrate the approach numerically on displacement interpolation and Wasserstein barycenters, utilizing a structure-preserving discretization of the Heisenberg group and convolutional entropic optimal transport.
\end{abstract}
 




\section{Introduction}

Time-frequency representations describe how the oscillatory content of a signal is distributed jointly in time and frequency. For $f,\psi\in L^2(\mathbb R;\mathbb C)$, we use the short-time Fourier transform
\begin{equation}
G_\psi f(t,\omega)=\int_{\mathbb R} f(\tau)\,\overline{\psi(\tau-t)}e^{-i\omega(\tau-t)}\,d\tau
=A_f^\psi(t,\omega)e^{i\Theta_f^\psi(t,\omega)}.
\end{equation}
Here $A_f^\psi\coloneq|G_\psi f|$ is the STFT magnitude; the spectrogram is $(A_f^\psi)^2$. Previous optimal transport-based approaches \cite{cazelles2020, valdivia2025} to time-frequency analysis act only on a nonnegative representation such as $A_f^\psi$ or $(A_f^\psi)^2$. This discards the phase $\Theta_f^\psi$, even though phase is the variable that reflects the coherent interplay of time shifts, frequency modulations, and global phase rotations.

The geometry that combines these three transformations is the three-dimensional Heisenberg group. Its Schr\"odinger representation realizes time-frequency shifts as a unitary group action, while its contact distribution expresses the classical phase-frequency consistency relation $d\phi=\omega\,dt$. This suggests replacing a planar time-frequency density by a probability measure on the reduced Heisenberg group, in which the central coordinate $\phi$ is periodic and corresponds to the STFT phase. Such a lift is equivariant under the Schr\"odinger action; consequently, every Wasserstein distance built from a left-invariant metric is invariant under common delay, Doppler and phase transformations.

This point of view is rooted in the classical role of the Heisenberg group in harmonic analysis \cite{Thangavelu1998, Howe1980, Koranyi1983, Folland89}. Time and frequency shifts do not commute: their commutator is a scalar phase rotation, and the Heisenberg group is precisely the central extension that turns this projective structure into a group representation. The resulting geometry is therefore not an auxiliary construction imposed on the time-frequency plane, but is intrinsic to the covariance properties of the Gabor transform itself. 
On the geometric side, the Heisenberg group is also the basic example of a noncommutative Carnot group and of a subRiemannian manifold \cite{Capogna, Ju2009}. Its horizontal distribution is bracket-generating, its natural metric is the Carnot-Carath\'eodory distance, and the anisotropic scaling between horizontal and central variables captures the fact that phase accumulates at second order with respect to horizontal motion.

At the differential level, this geometry appears through instantaneous frequency. If a complex oscillation is locally written as $z(t)=A(t)e^{i\phi(t)}$,
its instantaneous frequency is defined as $\hat\omega(t)=\dot \phi(t)$, so the lifted curve $\gamma(t)=\left(t, \hat\omega(t), \phi(t)\right)$
satisfies $(d\phi-\omega\,dt)(\dot{\gamma})=0$:
phase-frequency consistency is horizontality for the standard contact
structure.
For a general Gabor transform, the intersection of the phase surface with the horizontal distribution singles out a horizontal direction. Its projection onto the time-frequency plane is the displacement field of differential reassignment \cite{Cha97,daudet99}, the classical method that uses derivatives of the phase to relocate energy towards more concentrated representations, yielding a type of sharpening of the time-frequency representation which can be used, for instance, towards the extraction of instantaneous frequencies.

A closely related geometric viewpoint was developed in \cite{DUITS2013}, which emphasized that the reduced Heisenberg group is the natural domain of a Gabor transform when its phase covariance is retained. There, left-invariant evolution equations on the reduced Heisenberg group induce translation- and modulation-covariant processing of signals; in particular, nonlinear convection equations give
a geometric formulation of differential reassignment. We start from the same covariance principle, but the representation and the geometry are different. Rather than evolving the complex Gabor transform, we associate with a signal a probability measure supported on its phase graph and compare such measures by optimal transport. This combines, in a single representation, the contact geometry governing phase coherence and the metric geometry governing displacement of mass.
The lift also differs from other higher-dimensional lifts, such as time-frequency-chirp representations \cite{Chen23}, which add a chirp-rate coordinate to separate modes with crossing instantaneous frequencies. There the extra coordinate is an additional local descriptor; here it is the phase already present in the complex STFT, with a transformation law dictated by the Schr\"odinger representation.

Optimal transport provides the second ingredient of our construction. On Riemannian spaces, Wasserstein distances quantify discrepancies between probability measures by the geometric cost of moving mass along geodesics and have become useful tools for comparing, interpolating, and averaging structured distributions.
An extension to subRiemannian geometries has been studied in \cite{figalli2010, citti2025}, while for the Heisenberg group specifically in \cite{AMBROSIO2004, Rigot2005, deprig10, FIGALLI2008, pass2020}. In particular, \cite{AMBROSIO2004} proves existence and uniqueness of optimal transport maps for the quadratic cost when the source measure is absolutely continuous, and \cite{deprig10} treats Monge's
problem with the linear cost.
The geometry of Wasserstein space over the Heisenberg group differs markedly from the Riemannian case: for instance, the Heisenberg group does not satisfy any curvature-dimension condition \cite{Ju2009}.

Optimal transport has also begun to enter signal and time-frequency analysis directly. The Wasserstein-Fourier distance \cite{cazelles2020}, for example, compares stationary time series by transporting their normalized spectral energy distributions. More recently, Wasserstein transport and barycenters of spectrograms have been used for interpolation of audio signals \cite{valdivia2025}. These methods demonstrate the ability of transport to compare oscillatory content in a way that is sensitive to displacement in frequency or time-frequency space rather than merely pointwise variation. At the same time, when transport is performed on a power spectrum, magnitude STFT, or spectrogram, the transported object is necessarily nonnegative and the complex phase is absent. Phase must then either be ignored or reconstructed separately. Our aim is to build the phase into the geometry before the transport problem is posed.
Closest in spirit is \cite{pai2025}, where images are lifted to the group of roto-translations and transported with left-invariant costs, so that transport respects the group action. We develop the analogous principle for time-frequency analysis, with a substantially different geometric object: the symmetry group is the Heisenberg group, the lifted coordinate is the complex phase rather than an orientation, and the horizontal distribution encodes phase-frequency consistency. 
The noncommutative geometry thus not only provides an invariant transport cost; it also describes which infinitesimal variations of time, frequency
and phase are mutually coherent.

\textbf{Contributions:}
\begin{enumerate}
    \item \emph{A phase-equivariant lift.} We associate with a nonzero signal $f$ a probability measure $\mu^\psi_f$ on the reduced Heisenberg group $\mathring{\mathbb H}$ (Definition \ref{def: meas}). Its time-frequency marginal is a normalized STFT magnitude, with a cutoff near the zero set of $G_\psi f$, where the phase is unstable \cite{Balasz16}; its central coordinate is the phase of $G_\psi f$. The Schr\"odinger action on the signal becomes left translation, and the action on the window becomes right translation (Theorem \ref{thm:lift-equivariance}):
    \begin{equation}
        \mu^\psi_{U_hf}=(L_h)_\#\mu^\psi_f,\qquad \mu^\psi_f=(R_h)_\#\mu^{U_h\psi}_f .
    \end{equation}
    Hence, for any left-invariant distance on $\mathring{\mathbb H}$, such as the Carnot-Carath\'eodory or the Kor\'anyi distance, $W_p(\mu^\psi_{U_hf},\mu^\psi_{U_hg})=W_p(\mu^\psi_f,\mu^\psi_g)$ (Corollary \ref{cor:wasserstein-invariance}): two signals have the same Wasserstein discrepancy after a common delay, frequency modulation or phase rotation. This invariance follows from the ambient geometry rather than being imposed on the transport problem. Moreover, horizontal right translations are optimal transport maps (Corollary \ref{cor: monge}), and in radar terms an echo is the left translate of the lifted auto-ambiguity function by its delay, Doppler shift and phase (Remark \ref{rmk: radar}).
    \item \emph{Reconstruction and stability.} Although $\mu^\psi_f$ is a positive measure, a linear synthesis operator recovers $f$, up to its normalization, with a uniform error controlled by the low-amplitude coefficients discarded by the cutoff (Proposition \ref{prop: rec}). Without cutoff, $\mu^\psi_f$ and $\psi$ determine $f$ up to a positive factor (Corollary \ref{cor: rec1}), and so does synthesis with any window $\psi'$ such that $\langle\psi,\psi'\rangle\neq0$ (Lemma \ref{lemma: rec2}). Conversely, synthesis is Lipschitz from $W_1$ to $L^2$, and to $L^\infty$ for regular windows (Proposition \ref{prop:Rlip2}). The $L^2$ constant is minimized when the metric scale $\lambda$ matches the time-frequency spread of the window. Wasserstein closeness of lifts therefore implies closeness of the synthesized signals.
    \item  \emph{Reassignment as a Wasserstein gradient flow, and its non-convexity.} For Gaussian windows, differential reassignment lifts to the (formal) Wasserstein gradient flow of $\mathcal V_f(\mu)=-\int\log A_f\circ\mathrm{pr}\dd\mu$ on $\mathring{\mathbb H}$, with $\mathrm{pr}$ denoting the projection onto the first two coordinates (Proposition \ref{prop:gradient}). For every nonzero signal, $\mathcal V_f$ is not geodesically convex (Proposition \ref{prop: non-conv}).
    \item \emph{Computation of OT on the reduced Heisenberg group.} On compatible Heisenberg lattices the Gibbs kernel depends only on group difference, so the Sinkhorn kernel multiplications are restrictions of discrete group convolutions. We use this structure for transport interpolation and entropy-regularized barycenters.
\end{enumerate}







The remainder of the paper is organized as follows. Section \ref{sec: Hei} recalls the Heisenberg group, its reduced quotient and subRiemannian geometry, the Schr\"odinger representation and the covariance of the Gabor transform. Section \ref{sec: OT} reviews optimal transport on Lie groups with left-invariant distances, its entropic regularization, the resulting convolutional algorithms on $\mathring{\mathbb H}$, and the optimality of horizontal right translations. Section \ref{sec: main} introduces the phase lift and studies its equivariance, its reconstruction and Wasserstein stability, and the interpretation of differential reassignment as a Wasserstein gradient flow together with the failure of geodesic convexity. Section \ref{sec: num} presents numerical experiments on displacement interpolation and Wasserstein barycenters. Appendix \ref{app: sym} collects the symmetric model of the Heisenberg group, its geodesics, the Kor\'anyi gauge and the isometry $F_\lambda$ between the two models.

\section{The Heisenberg group}\label{sec: Hei}
In this section, we define the Heisenberg group in the coordinate system used throughout the paper, its reduced quotient, its subRiemannian geometry, and the Schr\"odinger representation relevant to time-frequency analysis.

Let $\Hb \coloneq \R^3_{(t,\omega,\phi)}$ and define the differential 1-form, and codimension 1-tangent distribution
\begin{equation}
  \alpha \coloneq d\phi - \omega\,dt, \qquad D \coloneq \ker\alpha.
\end{equation}

\begin{proposition}\label{prop:jet}
The following statements hold.
\begin{enumerate}[label=(\roman*)]
  \item $\alpha$ is a contact form:
        $\alpha \wedge d\alpha = d\phi \wedge dt \wedge d\omega \neq 0$.
        Its Reeb field is $R_\alpha \coloneq \partial_\phi$.
  \item The tangent vector fields
        $X_1 = \partial_t + \omega\,\partial_\phi$ and $X_2 = \partial_\omega$ constitute a global frame for $D$. Moreover $[X_1, X_2] = -\partial_\phi = -R_\alpha$, so $D$ has growth vector $(2,3)$, is bracket-generating of step~$2$, and is equiregular.
  \item\label{it: scale} For $\lambda>0$, declaring $X_1^\lambda\coloneq\lambda X_1, X_2^\lambda\coloneq\lambda^{-1}X_2$ orthonormal defines a one-parameter family of subRiemannian metrics $(\cdot, \cdot)_{\Hb,\lambda}$ on $\Hb$. 
  \item The binary operation
        \begin{equation}\label{eq:jet-group}
          (t,\omega,\phi) \star (t',\omega',\phi')
          \;=\;
          \bigl(t+t',\; \omega+\omega',\; \phi+\phi'+\omega t'\bigr)
        \end{equation}
        gives $\Hb$ a group structure and $X_1, X_2, R_\alpha,
        D, \alpha$ are all left-invariant. The inverse is 
        \begin{equation}
            (t,\omega, \phi)^{-1}=(-t,-\omega, -\phi +\omega t).
        \end{equation}
  \item The anisotropic dilations $\delta_r(t,\omega,\phi) = (rt, r\omega, r^2\phi)$
        for $r > 0$ are group automorphisms satisfying $\delta_r^*\alpha = r^2\alpha$.
\end{enumerate}
\end{proposition}

\begin{proof}
Since $d\alpha = -d\omega \wedge dt = dt \wedge d\omega$, the contact identity follows directly. It is also immediate that $R_\alpha$ satisfies the two defining properties of the Reeb field, namely $\alpha(R_\alpha)=1$ and $\iota_{R_\alpha}d\alpha = 0$.
Direct substitution gives $\alpha(X_1) = \alpha(X_2) = 0$ and $[X_1,X_2] = -\partial_\phi$.
For fixed $x = (t,\omega,\phi)$, the differential of
left translation $L_x(y) = x\star y$ sends
$\partial_{t'} \mapsto X_1$, $\partial_{\omega'} \mapsto X_2$, $\partial_{\phi'} \mapsto \partial_\phi$,
proving left-invariance. Associativity of \eqref{eq:jet-group} follows from the
bilinearity of the cocycle $(\omega, t') \mapsto \omega t'$. The dilation assertions
hold by direct substitution. 
\end{proof}
We denote by $d_{\Hb,\lambda}$ the subRiemannian distance in $\Hb$ associated with the metric $(\cdot, \cdot)_{\Hb,\lambda}$. Moreover, we denote by $L_h, R_h\colon  \Hb\to \Hb$ the products $L_h x\coloneq h\star x$, and $R_h x\coloneq x\star h$ for any $h\in \Hb$.

An alternative distance in the Heisenberg group with an explicit algebraic expression is the one derived from the Kor\'anyi gauge introduced in Appendix \ref{app: sym}, which for $x=(t,\omega,\phi)$ and $y=(t',\omega',\phi')$ results in
\begin{equation}\label{eq: koranyiintro}
d_{K,\lambda}(x,y) := \left[ \left( \frac{(t-t')^2}{\lambda^2} + \lambda^2 (\omega-\omega')^2 \right)^{2} + 16 \left( \phi - \phi' - \tfrac{1}{2}(t-t')(\omega+\omega') \right)^{2} \right]^{1/4}.
\end{equation}

For later use in the numerical treatment (Section \ref{subsec: computational}), we record which uniform grids are compatible with the group structure. Let $\Delta t,\Delta\omega,\Delta\phi>0$ and
\begin{equation}
  \calL_\Delta\coloneq\Delta t\,\mathbb{Z}\times\Delta\omega\,\mathbb{Z}\times\Delta\phi\,\mathbb{Z}\subset\mathbb{H}.
\end{equation}
Since 
\begin{equation}
    (k\Delta t,l\Delta\omega,m\Delta\phi)\star(k'\Delta t,l'\Delta\omega,m'\Delta\phi)=\bigl((k+k')\Delta t,(l+l')\Delta\omega,(m+m')\Delta\phi+lk'\,\Delta t\Delta\omega\bigr),
\end{equation}
the set $\calL_\Delta$ is closed under $\star$ if and only if $lk'\,\Delta t\Delta\omega\in\Delta\phi\,\mathbb{Z}$ for all $l,k'\in\mathbb{Z}$, and the case $l=k'=1$ is already necessary and sufficient. As $(t,\omega,\phi)^{-1}=(-t,-\omega,-\phi+\omega t)$ and $\omega t=kl\,\Delta t\Delta\omega$, inverses
are then automatic. Hence $\calL_\Delta$ is a subgroup of $\mathbb{H}$ if and only if
\begin{equation}\label{eq: structH}
    \frac{\Delta t\,\Delta\omega}{\Delta\phi}\in\mathbb{Z}.
\end{equation}

The Heisenberg group is sometimes defined through an alternative coordinate system. These coordinates are better adapted to standard formulas for subRiemannian geodesics and homogeneous gauges, but are less suitable for time-frequency analysis. In Appendix \ref{app: sym}, we recall this alternative model which we denote as $\bG$, review its basic properties, and pass between $\bG$ and $\bH$ by an explicit Carnot-group isomorphism.

\subsection{The reduced group}

For signal analysis the phase is defined modulo $2\pi$. We therefore quotient the center by the discrete subgroup
\begin{equation}
2\pi\Z_\phi\coloneq\{(0,0,2\pi k):k\in\mathbb Z\}.
\end{equation}
Since $2\pi\Z_\phi$ is a normal subgroup, 
\begin{equation}
    \mathring \Hb\coloneq\Hb/2\pi \Z_\phi
\end{equation}
is a Lie group. The action of $2\pi \Z_\phi$ on $\Hb$ is by deck transformations, each of which is a left translation and hence preserves the left-invariant contact and subRiemannian structures. Consequently $\alpha$, $D$, and the horizontal metric descend uniquely to $\mathring \Hb$, and the quotient map is a local subRiemannian isometry. The quotient distance can be written as
\begin{equation}
    \mathring d_{\Hb,\lambda}([x],[y])=\inf_{k\in\mathbb Z}d_{\Hb,\lambda}\left(x,y\star(0,0,2\pi k)\right).
\end{equation}
Similarly, the Kor\'anyi distance becomes
\begin{equation}\label{eq:koranyi_reduced}
    \mathring d_{K,\lambda}([x],[y])\coloneq\inf_{k\in\mathbb Z} d_{K,\lambda}\bigl(x,y\star(0,0,2\pi k)\bigr).
\end{equation}
The isomorphism $F_\lambda$ defined in Proposition \ref{prop:equi} preserves the central lattice and therefore descends to an isomorphism between the corresponding reduced $\mathring\bG$ and physical-phase $\mathring \Hb$ models.

\begin{proposition}\label{prop: wrapped geodesics}
    Let $[x]\neq[y]$ in $\mathring \Hb$ be represented by $x=(t,\omega, \phi)$ and $y=(t',\omega',\phi')$.
    \begin{enumerate}
        \item [(i)] $k^*\in \Z$ is optimal for $\mathring d_{\Hb,\lambda}([x],[y])$ if and only if $|\phi'-\phi-\frac{1}{2}(\omega+\omega')(t'-t)+2\pi k^*|\leq \pi$.
        \item [(ii)] The projection onto $\rH$ of a geodesic $\gamma$ between $x$ and $y\star (0,0,2\pi k^*)$ is a geodesic in $\rH$.
    \end{enumerate}
    \begin{proof}
        (i) Since $(0,0,2\pi k)$ is central and fixed by $F_\lambda$, left invariance gives
        $d_{\Hb,\lambda}\big(x,y\star(0,0,2\pi k)\big)=d_\bG\big(0,((t'-t)/\lambda,\lambda(\omega'-\omega),\delta+2\pi k)\big)$, where $\delta=\phi'-\phi-\frac12(\omega+\omega')(t'-t)$. By Proposition \ref{prop:geodesics}, this distance increases strictly with $|\delta+2\pi k|$, so it is minimal if and only if $|\delta+2\pi k|\le\pi$. \\ 
        (ii) The projection $\pi\colon \Hb\to\rH$ is $1$-Lipschitz, so by (i) $\mathring d_{\Hb,\lambda}([x],[y])\le \ell(\pi\circ\gamma)\le \ell(\gamma)=d_{\Hb,\lambda}\big(x,y\star(0,0,2\pi k^*)\big)=\mathring d_{\Hb,\lambda}([x],[y])$.   
    \end{proof}
\end{proposition}
\begin{remark}
    Since the optimal phase defect satisfies $|\phi'-\phi-\frac12(\omega+\omega')(t'-t)+2\pi k^*|\leq\pi$, the area enclosed by a corresponding geodesic in $\bG$ in the $(t,\omega)$-plane is at most $\pi$. Geodesics of $\rH$ therefore follow straighter lines the farther apart $x$ and $y$ are in time and frequency: combining the area bound with the explicit geodesic formula in Proposition \ref{prop:geodesics}, each geodesic has a turning angle of at most $12\pi/r^2$ and $\mathring d_{\Hb,\lambda}([x],[y])^2=r^2+O(r^{-2})$, with $r=|(\lambda^{-1}(t'-t),\lambda(\omega'-\omega))|$. Thus, the phase difference plays a more relevant role the closer the points are in the $(t,\omega)$-domain.
\end{remark}

Regarding the structure-preserving discretization in the reduced group, the condition in Equation \eqref{eq: structH}  has the additional constraint $\Delta\phi=2\pi/N_\phi$ with $N_\phi\in\mathbb{N}$. Given this, the lattice becomes a discrete subgroup of $\mathring{\mathbb{H}}$ if and only if
\begin{equation}\label{eq: structHred}
    \frac{N_\phi \Delta t \Delta\omega}{2\pi}\in \Z.
\end{equation}
We will also assume that the grid aspect ratio $\Delta t/\Delta\omega$ is consistent with the scaling $\lambda$, i.e. $\Delta t/\Delta\omega\approx \lambda^2$.

\subsection{Application to signal analysis}
The Heisenberg group enters time-frequency analysis through the basic noncommutativity of time translations and frequency modulations. Indeed, with the conventions used below in Definition \ref{def:tf-shifts}, 
\begin{equation} 
T_s M_\eta = e^{-i\eta s} M_\eta T_s . 
\end{equation} 
Thus time and frequency shifts form only a projective representation of $\mathbb R^2$: changing their order produces a global phase factor. The Heisenberg group resolves precisely this obstruction by adjoining phase as a central variable. Its Schr\"odinger representation then realizes translations, modulations, and global phase rotations as a single unitary group action.

There is a complementary differential-geometric interpretation. If a complex oscillation is locally written as $A(t)e^{i\phi(t)}$, then its instantaneous frequency is the derivative of its phase. The relation 
\begin{equation} 
d\phi = \omega\,dt 
\end{equation} 
is exactly the horizontality condition associated with the contact form $\alpha=d\phi-\omega\,dt$ introduced above. Hence the algebraic structure of time-frequency shifts and the differential relation between phase and instantaneous frequency are encoded by the same Heisenberg geometry. We make these two observations precise in the following subsections. They will provide the geometric basis for the phase lift of the Gabor transform introduced later.

We endow the time-frequency plane $(t,\omega)$ with the scaled Euclidean metric $g_\lambda=\lambda^{-2}  dt^2+\lambda^2 d\omega^2$, and we denote by $|z|_{g_\lambda}^2=z_1^2/\lambda^2+\lambda^2z_2^2$ the corresponding norm. This immediately yields the following lemma, which also holds for the Kor\'anyi metric.
\begin{lemma}\label{lemma: projection}
    The plane projection $\operatorname{pr}\colon (\rH, \mathring d_{\Hb,\lambda})\to (\R^2, |\cdot|_{g_\lambda})$ onto the first two coordinates is $1$-Lipschitz.
\end{lemma}
\begin{remark}
    The parameter $\lambda>0$ appearing in the subRiemannian metric and in the time-frequency metric reflects the fact that $t$ has units of time and the angular frequency $\omega$ has reciprocal units. So the rescaled coordinates $t/\lambda$ and $\lambda \omega$ are dimensionless, and $g_\lambda$ measures time and frequency displacements relative to the temporal resolution $\lambda$ and the reciprocal resolution $1/\lambda$ in frequency. The parameter will be later chosen consistently with the scale of the Gaussian window in the Gabor transform.
\end{remark}


\subsubsection{Time-frequency shifts and the Schr\"odinger representation}

We first recall the elementary operators generating time-frequency shifts on $L^2(\mathbb R;\mathbb C)$.

\begin{definition}
\label{def:tf-shifts}
For $(s,\eta)\in\mathbb R^2$, we define the time translation and frequency modulation operators by 
\begin{equation} 
(T_s f)(\tau) \coloneq f(\tau-s), \qquad (M_\eta f)(\tau) \coloneq e^{i\eta\tau}f(\tau), 
\end{equation} 
and the time-frequency shift as
\begin{equation} 
\Pi_{s,\eta}\coloneq T_s\circ M_\eta . 
\end{equation}
\end{definition}

\begin{definition}
\label{def:phase-rotation}
For $\theta\in\mathbb R/2\pi\mathbb Z$, we define the global phase rotation 
\begin{equation} 
P_\theta f\coloneq e^{i\theta}f . 
\end{equation}
\end{definition}
The operators $T_s$ and $M_\eta$ are unitary but do not commute. More precisely, 
\begin{equation} \label{eq:weyl-commutation} 
T_sM_\eta = e^{-i\eta s}M_\eta T_s . 
\end{equation} 
The scalar factor in \eqref{eq:weyl-commutation} is exactly the central coordinate appearing in the Heisenberg group law. This is the classical motivation for the Schr\"odinger representation; see \cite[Section 1.3]{Folland89} or \cite[Section~1.2]{Thangavelu1998} for a comprehensive account.

\begin{definition} \label{def:schrodinger} 
For $h=(s,\eta,\theta)\in\mathring \Hb$, define 
\begin{equation} \label{eq:schrodinger} 
U_h f(\tau) \coloneq e^{i\theta}e^{i\eta(\tau-s)}f(\tau-s) = (P_\theta\circ T_s\circ M_\eta)f(\tau). 
\end{equation} 
\end{definition}
The map 
\begin{equation} 
U\colon\mathring \Hb\longrightarrow \mathcal U\bigl(L^2(\mathbb R;\mathbb C)\bigr), \qquad h\longmapsto U_h , 
\end{equation} 
is a unitary representation of the reduced Heisenberg group. In particular, for all $x,y\in\rH$, 
\begin{equation} \label{eq:schrodinger-group-action} 
U_xU_y=U_{x\star y}, \qquad U_{h^{-1}}=U_h^{-1}, 
\end{equation}
and 
\begin{equation} \label{eq:schrodinger-unitary} 
\|U_hf\|_{L^2(\mathbb R;\mathbb C)} = \|f\|_{L^2(\mathbb R;\mathbb C)}. 
\end{equation} 
The representation identity follows directly from the group law. Indeed, if $x=(s,\eta,\theta)$ and $y=(s',\eta',\theta')$, then 
\begin{align} 
U_xU_yf(\tau) &= e^{i(\theta+\theta')} e^{i\eta(\tau-s)} e^{i\eta'(\tau-s-s')} f(\tau-s-s') \\ 
&=U_{x\star y}f(\tau). 
\end{align} 
Unitarity follows from the fact that translations, modulations, and phase rotations are individually unitary.

The infinitesimal form of this representation identifies the horizontal directions of $\Hb$ with the standard differentiation and multiplication operators of time-frequency analysis.
\begin{lemma} \label{lemma:U-deriv} 
Let 
\begin{equation} 
\mathcal U_f\colon\mathring \Hb\longrightarrow L^2(\mathbb R;\mathbb C), \qquad \mathcal U_f(h)\coloneq U_hf 
\end{equation} 
be the orbit map of $f$. Assume 
\begin{equation} 
f\in W^{1,2}(\mathbb R;\mathbb C), \qquad \tau f(\tau)\in L^2(\mathbb R;\mathbb C). 
\end{equation} 
Then, for every $h\in\rH$, 
\begin{equation} \label{eq:infinitesimal-schrodinger} 
(X_1\mathcal U_f)(h) =U_h(-\partial_\tau f), \qquad (X_2\mathcal U_f)(h) =U_h(i\tau f), \qquad (R_\alpha\mathcal U_f)(h) =U_h(if). 
\end{equation} 
Equivalently, the derived representation at the identity is 
\begin{equation} \label{eq:derived-representation} 
dU(X_1)=-\partial_\tau, \qquad dU(X_2)=i\tau, \qquad dU(R_\alpha)=i\,\mathrm{Id}. 
\end{equation} 
\end{lemma}
\begin{proof}
Since $X_1,X_2,R_\alpha$ are left-invariant, their action on the orbit map can be computed by differentiating curves of the form $h\star\exp(\epsilon X)$.
        For $X_1$, 
        \begin{align} 
        (X_1\mathcal U_f)(h) &= \left. \frac{d}{d\epsilon} \right|_{\epsilon=0} U_{h\star(\epsilon,0,0)}f \\ 
        &= U_h \left( \left. \frac{d}{d\epsilon} \right|_{\epsilon=0} f(\tau-\epsilon) \right) = U_h(-\partial_\tau f). 
        \end{align}
        The same technique can be used for the other two vector fields.
    \end{proof}

\subsubsection{Horizontality is phase-frequency consistency}
We next relate the contact structure of $\Hb$ to the classical notion of instantaneous frequency. For a real-valued oscillatory signal, a standard construction first forms the associated analytic signal using the Hilbert transform and then defines the instantaneous frequency as the derivative of its phase; see, for instance, \cite{Huang2009OnIF} and \cite[Section~1]{Thakur_2011}. In the present work the signals are already complex-valued, so no analytic extension is required: whenever 
\begin{equation} 
z(t)=A(t)e^{i\phi(t)}, \qquad A(t)>0, 
\end{equation} 
the local instantaneous frequency is simply $\phi'(t)$. 
The Heisenberg contact structure gives a geometric formulation of precisely this relation.

\begin{proposition}
\label{prop:horizontal}
Let $\gamma(s) = (t(s), \omega(s), \phi(s))$ be an absolutely continuous curve in $\Hb$, with $\phi$ an absolutely continuous lift of its phase coordinate. Then
\begin{equation}\label{eq:horizontal-phase}
  \gamma \text{ is } D\text{-horizontal}
  \;\;\Longleftrightarrow\;\;
  \dot\phi(s) = \omega(s)\,\dot t(s) \quad \text{for a.e.~} s.
\end{equation}
In particular, if the curve is parametrized by physical time, 
\begin{equation} 
\gamma(t)=\bigl(t,\omega(t),\phi(t)\bigr), 
\end{equation} 
then horizontality is equivalent to 
\begin{equation} \label{eq:if-horizontal} 
\omega(t)=\dot\phi(t). 
\end{equation} 
Consequently, for any $A\in C^1$ with $A(t)>0$ and 
\begin{equation} 
z(t)=A(t)e^{i\phi(t)}, 
\end{equation} 
one has 
\begin{equation} \label{eq:logarithmic-if} 
\operatorname{Im}\frac{\dot z(t)}{z(t)} = \dot \phi(t) = \omega(t). 
\end{equation} 
\end{proposition}

\begin{proof}
By definition, $\gamma$ is horizontal if and only if $\alpha(\dot\gamma)=0$. Since 
$\alpha=d\phi-\omega\,dt,$ we have 
\begin{equation} 
\alpha(\dot\gamma) = \dot\phi-\omega\dot t, 
\end{equation} 
which proves 
\eqref{eq:horizontal-phase}. 
If $t$ itself parametrizes the curve, then $\dot t=1$, giving 
\eqref{eq:if-horizontal}. 
For $z=Ae^{i\phi}$ with $A>0$, 
\begin{equation} \frac{\dot z}{z} = \frac{\dot A}{A}+i\dot \phi, 
\end{equation} and taking the imaginary part gives \eqref{eq:logarithmic-if}.
\end{proof}

Proposition \ref{prop:horizontal} is the geometric principle that will be used throughout the signal-analysis constructions below. A phase-valued time-frequency representation naturally defines a surface in the reduced Heisenberg group. That surface need not itself be horizontal; rather, its intersection with the contact distribution selects directions along which variations of time, frequency, and phase are mutually consistent. For the Gabor transform, these distinguished horizontal directions will be shown to recover the displacement field of differential reassignment, thereby linking the classical phase-derivative construction to the intrinsic geometry of $\Hb$.

\subsubsection{Covariance and regularity of Gabor transform}
We now turn to Gabor transforms of $L^2(\R, \C)$ signals, defined using the corresponding inner products $\langle\cdot,\cdot\rangle$ which are linear in the first and conjugate-linear in the second argument.

\begin{definition}\label{def:gabor}
Let $\psi\in L^2(\R,\bC)\setminus\{0\}$. The short-time Fourier (or Gabor) transform of $f\in L^2(\R,\bC)$ is
\begin{equation}\label{eq:gabor}
G_\psi f(t,\omega)\coloneq\langle f,U_{(t,\omega,0)}\psi\rangle
=\int_\R f(\tau)\,\overline{\psi(\tau-t)}\,e^{-i\omega(\tau-t)}\,d\tau
=A_f^\psi(t,\omega)\,e^{i\Theta_f^\psi(t,\omega)},
\end{equation}
with amplitude $A_f^\psi\coloneq |G_\psi f|$ and phase $\Theta_f^\psi\in\R/2\pi\Z$ defined on $\{A_f^\psi>0\}$. The Heisenberg-Gabor transform (cf. \cite{DUITS2013}) is the function on $\mathring\Hb$
\begin{equation}\label{eq: full-gabor}
W_\psi f(t,\omega,\phi)\coloneq\langle f,U_{(t,\omega,\phi)}\psi\rangle=e^{-i\phi}\,G_\psi f(t,\omega),
\end{equation}
where the second equality uses $U_{(t,\omega,\phi)}=e^{i\phi}U_{(t,\omega,0)}$.
\end{definition}

\begin{lemma}\label{lem:Wcov}
For all $f,\psi\in L^2(\R,\bC)$ and $h\in\mathring\Hb$,
\begin{equation}\label{eq:Wcov}
W_\psi(U_hf)=W_\psi f\circ L_h^{-1},
\qquad
W_{U_h\psi}f=W_\psi f\circ R_h.
\end{equation}
Equivalently, for $h=(s,\eta,\theta)$,
\begin{align}\label{eq:Gcov}
&G_\psi(U_hf)(t,\omega)=e^{i(\theta+\eta(t-s))}\,G_\psi f(t-s,\omega-\eta),
\\
&G_{U_h\psi}f(t,\omega)=e^{-i(\theta+\omega s)}\,G_\psi f(t+s,\omega+\eta).
\end{align}
\end{lemma}
\begin{proof}
By unitarity and \eqref{eq:schrodinger-group-action}, $U_x^*=U_{x^{-1}}$, so for every $y\in\mathring\Hb$
\begin{equation}
W_\psi(U_xf)(y)=\langle f,U_{x^{-1}}U_{y}\psi\rangle=W_\psi f(x^{-1}\star y),
\qquad
W_{U_x\psi}f(y)=\langle f,U_{y}U_x\psi\rangle=W_\psi f(y\star x).
\end{equation}
For $y=(t,\omega,0)$ one has $x^{-1}\star y=(t-s,\omega-\eta,-\theta-\eta(t-s))$ and $y\star x=(t+s,\omega+\eta,\theta+\omega s)$, and \eqref{eq: full-gabor} gives the desired relations.
\end{proof}

Thus the Schr\"odinger action on the signal (resp. window) acts on $W_\psi f$ by left (resp. right) translation. On $G_\psi f$ it translates the amplitude and shifts the phase by an affine function of $(t,\omega)$, which is exactly what the third coordinate absorbs.

\begin{lemma}\label{lemma: AC0}
Let $f,\psi\in L^2(\R,\bC)$.\\
\textup{(i)} $G_\psi f\in C_0(\R^2)$ and $\|G_\psi f\|_\infty\le\|f\|_2\|\psi\|_2$.\\
\textup{(ii)} For all $g,\psi'\in L^2(\R,\bC)$ (Moyal's identity),
\begin{equation}\label{eq:moyal}
\langle G_\psi f,G_{\psi'}g\rangle_{L^2(\R^2)}=2\pi\,\langle f,g\rangle\,\overline{\langle\psi,\psi'\rangle}.
\end{equation}
In particular, $\|G_\psi f\|_{L^2(\R^2)}^2=2\pi\|f\|_2^2\|\psi\|_2^2$ and, for $\psi\neq0$, weakly in $L^2(\R,\bC)$,
\begin{equation}\label{eq:gabor-inversion}
f=\frac{1}{2\pi\|\psi\|_2^2}\int_{\R^2}G_\psi f(t,\omega)\,U_{(t,\omega,0)}\psi\dd t\dd\omega.
\end{equation}
\textup{(iii)} For every $a>0$, $B_a\coloneq\{A_f^\psi\ge a\}$ is compact and $|B_a|\le (2\pi/a^2)\|f\|_2^2\|\psi\|_2^2$.
\end{lemma}
\begin{proof}
(i) The bound follows from Cauchy-Schwarz and unitarity of $U_*$ in the inner product. For $f,\psi\in C_c^\infty(\R)$, $G_\psi f$ is continuous and vanishes for $|t|$ large. Integrating by parts in $\tau$ gives $|G_\psi f(t,\omega)|\le C|\omega|^{-1}$ uniformly in $t$, so $G_\psi f\in C_0(\R^2)$. By sesquilinearity and the bound, $\|G_\psi f-G_{\psi'}f'\|_\infty\le\|f-f'\|_2\|\psi\|_2+\|f'\|_2\|\psi-\psi'\|_2$. The general case follows by density, since $C_0(\R^2)$ is closed under uniform limits.
(ii) Identity \eqref{eq:moyal} is \cite[Thm. 3.2.1]{Grochenig2001}; the factor $2\pi$ comes from the angular variable. Formula \eqref{eq:gabor-inversion} follows by pairing with an arbitrary $g \in L^2(\R,\bC)$, since $\langle U_{(t,\omega,0)}\psi,g\rangle=\overline{G_\psi g(t,\omega)}$.
(iii) $B_a$ is closed by continuity and bounded because $A_f^\psi$ vanishes at infinity. Chebyshev's inequality and (ii) give $a^2|B_a|\le\|A_f^\psi\|_{L^2(\R^2)}^2$.
\end{proof}

\section{Optimal transport on Lie groups}\label{sec: OT}

Optimal transport (OT) provides a principled and versatile framework for comparing measures defined on very general spaces, well beyond the Euclidean setting in which it originated: measures on Riemannian manifolds \cite{mccann2001polar}, on general metric-measure spaces \cite{ambrosio2005gradient}, and, closer to the present setting, on Lie groups carrying a left-invariant (sub-)Riemannian structure \cite{AMBROSIO2004,deprig10,figalli2010,Rigot2005}. This generality is one of the reasons OT has developed such rich and diverse connections throughout mathematics and its applications, from gradient flows and partial differential equations \cite{ambrosio2005gradient,Benamou1999computational} to geometric analysis \cite{lott2006some}, and, as we develop in this paper, to phase-equivariant time-frequency analysis.

This section focuses on the formulation and computation of OT problems for measures defined on a Lie group $G$ equipped with a left-invariant distance, and specializes them to the reduced Heisenberg group, with the Carnot-Carath\'eodory and Kor\'anyi distances $d = \mathring d_{\mathbb H,\lambda}$ and $d=\mathring d_{K,\lambda}$. Section~\ref{subsec: basic-defs} recalls the basic definitions and
properties of the Wasserstein distance on such a space; Section~\ref{subsec:entropic-reg} introduces its entropic regularization; Section~\ref{subsec: computational} discusses the resulting computational aspects and algorithms in the reduced Heisenberg group, including the reduction to group convolution and the Wasserstein-barycenter problem; and Section \ref{sec: right} focuses on instances where the right translation is an optimal transport map.

\subsection{Definitions and Properties}\label{subsec: basic-defs}
In this subsection, $G$ denotes a Lie group equipped with a
left-invariant distance $d$, i.e.\ $d(h\star x,h\star y) = d(x,y)$ for all
$h,x,y \in G$, where $\star$ denotes the group operation. Denote by $e\in G$ the neutral element.

\begin{definition}
\label{def:wasserstein}
Let $(G,d)$ be as above and let $p \geq 1$. For $\mu,\nu \in \mathcal
P(G)$, the $p$-Wasserstein distance between $\mu$ and $\nu$ is
\begin{equation}
W_p^p(\mu,\nu) \;\coloneq\; \inf_{\pi \in \Pi(\mu,\nu)} \int_{G \times G} d(x,y)^p \, \dd\pi(x,y),
\end{equation}
where $\Pi(\mu,\nu)$ is the set of couplings of $\mu$ and $\nu$, i.e.\
probability measures on $G \times G$ with first and second marginals
$\mu$ and $\nu$ respectively. The space of probability measures with finite $p$-moments (sometimes called the $p$-Wasserstein space) over $G$ is denoted as $\mathcal P_p(G) \coloneq \{\mu \in \mathcal P(G) : W_p(\mu,\delta_e) < \infty\}$.
\end{definition}

\begin{definition}
\label{def:barycenter}
Let $(G,d)$ and $p \geq 1$ be as above, let $n \geq 2$, and let
$\boldsymbol{\mu}^n = (\mu_1,\dots,\mu_n) \in \mathcal P_p(G)^n$. For a weight vector
$\mathbf{w}^n = (w_i)_{i=1}^n$ with $0 \leq w_i \leq 1$ and $\sum_{i=1}^n w_i = 1$, a $p$-Wasserstein barycenter of $\boldsymbol{\mu}^n$ with weights $\mathbf{w}^n$
is any minimizer
\begin{equation}
\mathrm{Bar}(\boldsymbol{\mu}^n,\mathbf{w}^n) \;\in\; \operatorname*{arg\,min}_{\nu \in \mathcal P_p(G)} \; \sum_{i=1}^n w_i \, W_p^p(\nu,\mu_i).
\end{equation}
\end{definition}

Theoretical properties of the barycenter problem with $G = \mathbb H$ are developed in \cite{pass2020}, while \cite{pai2025} deals with $G = SE(2)$. In our setting, existence of a barycenter can be proved by a direct application of the Prokhorov theorem, since both $(\rH, \mathring d_{\Hb, \lambda})$ and $(\rH, \mathring d_{K, \lambda})$ are proper metric spaces (i.e. in which every closed bounded subset is compact). Indeed the projection onto the first two coordinates is $1$-Lipschitz, while the phase coordinate belongs to the circle, which is itself compact. Moreover, by considering convex combinations of transport plans for each term in the sum, it is easy to see that this minimization problem is convex in $\nu$.

\begin{proposition}
\label{prop:Wp-left-invariant}
Let $(G,d)$ be a Lie group with left-invariant distance $d$, let $p \geq 1$,
and let $\mu,\nu \in \mathcal P_p(G)$. Then for every $h \in G$,
\begin{equation}
W_p\big((L_h)_\# \mu, \, (L_h)_\# \nu\big) \;=\; W_p(\mu,\nu).
\end{equation}
\end{proposition}

\begin{proof}
This equivariance of $W_p$ under left translation was established, in the
same generality, for the roto-translation group $\mathrm{SE}(2)$ in
\cite{pai2025}; the argument there extends verbatim to any Lie group
carrying a left-invariant distance, since it relies only on
left-invariance of $d$ and on the bijection between couplings of
$\mu,\nu$ and of $(L_h)_\#\mu,(L_h)_\#\nu$ induced by $L_h \times L_h$.
\end{proof}


\begin{proposition}
\label{prop:barycenter-equivariant}
Let $(G,d)$ be a Lie group with left-invariant distance $d$, let $p \geq
1$, let $n \geq 2$, let $\boldsymbol{\mu}^n = (\mu_1,\dots,\mu_n) \in \mathcal
P_p(G)^n$ and $\mathbf{w}^n$ as in Definition \ref{def:barycenter}. For every $h \in G$,
\begin{equation}
\mathrm{Bar}\big((L_h)_\# \boldsymbol{\mu}^n,\mathbf{w}^n\big) \;=\; (L_h)_\# \, \mathrm{Bar}(\boldsymbol{\mu}^n,\mathbf{w}^n),
\end{equation}
where $(L_h)_\#\boldsymbol{\mu}^n \coloneq \big((L_h)_\#\mu_1,\dots,(L_h)_\#\mu_n\big)$.
\end{proposition}

\begin{proof}
This is an immediate consequence of Proposition~\ref{prop:Wp-left-invariant}:
for any $\nu \in \mathcal P_p(G)$,
\begin{equation}
\sum_{i=1}^n w_i\, W_p^p\big((L_h)_\#\nu,(L_h)_\#\mu_i\big) \;=\; \sum_{i=1}^n w_i\, W_p^p(\nu,\mu_i),
\end{equation}
so $\nu \mapsto (L_h)_\#\nu$ is a bijection between the minimizers of the
two problems in Definition~\ref{def:barycenter}. The same equivariance
of barycenters under left translation was established for
$\mathrm{SE}(2)$ in \cite{pai2025}.
\end{proof}

\subsection{Entropic Regularization}\label{subsec:entropic-reg}

Realizing the generality of optimal transport computationally is a significant challenge, since the exact Kantorovich problem is a linear program whose size grows quadratically with respect to the number of points in the discretization. Many numerical strategies have been proposed to address this, the most popular being entropic regularization \cite{cuturi2013sinkhorn,benamou2015iterative,solomon2015convolutional}: adding a small entropy penalty makes the Kantorovich problem strictly convex, and the resulting problem can be solved by simple, efficient fixed-point iterations based on iterative projections to satisfy each marginal constraint. These iterations reduce to repeated applications of a single Gibbs kernel $K(x,y) = \exp(-d(x,y)^p/\varepsilon)$, where $d$ is the distance on the underlying space, and this strategy extends well beyond Euclidean domains. On Riemannian manifolds, and more generally on discretized geometric domains such as triangulated meshes and point clouds, the exact Gibbs kernel is rarely available in closed form and is typically replaced by the heat kernel of the Laplace-Beltrami operator or a discretization thereof to allow the same alternating scheme to be used \cite{solomon2015convolutional}. A related strategy on the Lie group of roto-translations instead replaces the group distance itself with a closed-form logarithmic distance approximation \cite{pai2025,pai2025entropic}. Moreover, once the kernel to be used is determined, such an alternating projection scheme also allows for the computation of Wasserstein barycenters of several measures \cite{benamou2015iterative}.

Fix $p\ge 1$ and $\varepsilon>0$. For $\mu,\nu\in\mathcal P_p(G)$, entropic regularization replaces the linear transport cost of the Kantorovich problem with a strictly convex one, penalizing the divergence of a coupling from the product measure $\mu\otimes\nu$:
\begin{equation}
W_{p,\varepsilon}^p(\mu,\nu) \;\coloneq\; \min_{\pi\in\Pi(\mu,\nu)} \int_{G\times G} d(x,y)^p\, d\pi(x,y) \;+\; \varepsilon\, \mathrm{KL}\big(\pi \,\big|\, \mu\otimes\nu\big),
\end{equation}
where $\mathrm{KL}(\pi\mid\mu\otimes\nu)\coloneq\int_{G\times G}\log\!\big(\tfrac{d\pi}{d(\mu\otimes\nu)}\big)\,d\pi$ when $\pi\ll\mu\otimes\nu$, and $+\infty$ otherwise. Unlike the unregularized problem, whose minimizers can concentrate on lower-dimensional subsets of $G\times G$, the problem above is strictly convex in $\pi$ and admits a unique minimizer $\pi_\varepsilon$ for every $\varepsilon>0$, converging to an optimal coupling of the original Kantorovich problem as $\varepsilon\to0$.

The same regularization extends directly to the barycenter problem of Definition~\ref{def:barycenter}. Given $\boldsymbol{\mu}^n=(\mu_1,\ldots,\mu_n)\in\mathcal P_p(G)^n$ and weights $\mathbf{w}^n=(w_1,\ldots,w_n)$, the entropic barycenter replaces each term $W_p^p(\nu,\mu_i)$ by its entropic counterpart, resulting in
\begin{equation}
\mathrm{Bar}_\varepsilon(\boldsymbol{\mu}^n,\mathbf{w}^n) \;\in\; \operatorname*{arg\,min}_{\nu\in\mathcal P_p(G)} \sum_{i=1}^n w_i\, W_{p,\varepsilon}^p(\nu,\mu_i).
\end{equation}
As in the two-marginal case, this turns a problem that is jointly convex in $\nu$ and the transport plans $(\pi_1,\ldots,\pi_n)$ for each term into one that is strictly convex, with a unique minimizer computable by the iterative Bregman projection scheme of \cite{benamou2015iterative}. This scheme alternates Kullback-Leibler projections of the couplings $\pi_1,\ldots,\pi_n$ onto their marginal constraints. For a single pair $(\mu,\nu)$, the minimizer $\pi_\varepsilon$ factorizes as $d\pi_\varepsilon(x,y) = u(x)\,K(x,y)\,v(y)\,d\mu(x)\,d\nu(y)$, where $K(x,y)=e^{-d(x,y)^p/\varepsilon}$ is the Gibbs kernel associated with the ground cost, and the scaling functions $u,v$ are obtained, up to a multiplicative constant, from the fixed-point recursion
\begin{equation}
u^{(k+1)}(x) = \left(\int_G K(x,y)\, v^{(k)}(y)\, d\nu(y)\right)^{-1}, \qquad
v^{(k+1)}(y) = \left(\int_G K(x,y)\, u^{(k+1)}(x)\, d\mu(x)\right)^{-1},
\end{equation}
which is also well known as the Sinkhorn algorithm \cite{cuturi2013sinkhorn}. Each iteration therefore costs two applications of the integral operator with kernel $K$, and the barycenter scheme of \cite{benamou2015iterative} reduces likewise to repeated applications of the same kernel across all $n$ marginals at once.

The efficiency of such methods relies on being able to compute these rescaling steps efficiently, and this is where the group structure of $G$ becomes computationally decisive: indeed, for the Euclidean distance one can just use FFT methods. In our case, since $d$ is left-invariant, $d(x,y)=d(e,x^{-1}\star y)$ for all $x,y\in G$, so the Gibbs kernel itself depends only on $x^{-1}\star y$: writing $k(z)=e^{-d(e,z)^p/\varepsilon}$, one has $K(x,y)=k(x^{-1}\star y)$, and each Sinkhorn update becomes a group convolution,
\begin{equation}
\int_G K(x,y)\, v(y)\, d\nu(y) = \int_G k(x^{-1}\star y)\, v(y)\, d\nu(y) = \big(k * (v\,\nu)\big)(x).
\end{equation}
Thus, any method that allows computing discrete group convolutions in an efficient manner yields an immediate benefit for obtaining the corresponding transport plans and barycenters.

\begin{remark}\label{rmk: entropic-equivariance}
Since $d$ is left-invariant and the Kullback-Leibler divergence is invariant under pushforward by any measurable bijection, in particular by the left translation $L_h\colon x\mapsto h\star x$ for $h\in G$, both terms in the definition of $W_{p,\varepsilon}^p$ transform consistently under a simultaneous left translation of $\mu,\nu$ and of the coupling $\pi$. Consequently
\begin{equation}
W_{p,\varepsilon}^p((L_h)_\#\mu, (L_h)_\#\nu) = W_{p,\varepsilon}^p(\mu,\nu), \qquad h\in G,
\end{equation}
as for the unregularized distance, and the same argument applied termwise to the entropic barycenter functional gives the equivariance
\begin{equation}
\mathrm{Bar}_\varepsilon\big(((L_h)_\#\mu_1,\ldots,(L_h)_\#\mu_n),\mathbf{w}^n\big) = (L_h)_\#\,\mathrm{Bar}_\varepsilon(\boldsymbol{\mu}^n,\mathbf{w}^n).
\end{equation}
Entropic regularization and the left-invariant geometry of $G$ are therefore compatible throughout.
\end{remark}

\subsection{Discretization and Algorithms}\label{subsec: computational}
Discrete measures on $\mathring{\mathbb H}$ take the form
\begin{equation}\label{eq: disc_measures}
\mu = \sum_{j=1}^J a_j\,\delta_{x_j}, \qquad \nu = \sum_{\ell=1}^L b_\ell\,\delta_{y_\ell}, \qquad x_j,y_\ell\in\mathring{\mathbb H},\ \ a,b \text{ probability vectors},
\end{equation}
and the entropic scheme of Section~\ref{subsec:entropic-reg} specializes to this setting in two different ways, according to whether the locations $x_j,y_\ell$ in the transport plans and barycenters are fixed in advance or are allowed to change in the computation.

\paragraph{Fixed support.} Here all measures under consideration share a common set of atoms, drawn from a fixed grid $\{z_1,\ldots,z_J\}\subset\mathring{\mathbb H}$, most commonly a uniform grid in the $(t,\omega)$ plane combined with a uniform sampling of the third coordinate making the resulting triples belong to a common discrete subgroup of $\mathring{\mathbb H}$ by satisfying the condition \eqref{eq: structHred}. In the fixed support regime only the weights vary, so $\mu=\sum_j a_j\delta_{z_j}$ and $\nu=\sum_\ell b_\ell\delta_{z_\ell}$. The discrete Gibbs kernel $K_{j\ell}=e^{-\mathring d_{\mathbb H,\lambda}(z_j,z_\ell)^p/\varepsilon}$ then inherits the left-invariance of $\mathring d_{\mathbb H,\lambda}$: since $\mathring d_{\mathbb H,\lambda}(z_j,z_\ell)=\mathring d_{\mathbb H,\lambda}(e,z_j^{-1}\star z_\ell)$, one has $K_{j\ell}=\kappa(z_j^{-1}\star z_\ell)$ for the single function $\kappa(h) = e^{-\mathring d_{\mathbb H,\lambda}(e,h)^p/\varepsilon}$ of the grid offset, and, since $\kappa(h)=\kappa(h^{-1})$, $K$ is symmetric. Every matrix-vector product appearing in the Sinkhorn recursion,
\begin{equation}
u_j^{(t+1)} = \frac{a_j}{(Kv^{(t)})_j}, \qquad v_\ell^{(t+1)} = \frac{b_\ell}{(Ku^{(t+1)})_\ell}, \qquad (Kv)_j=\sum_{\ell=1}^L \kappa(z_j^{-1}\star z_\ell)\,v_\ell = (\kappa * v)(z_j),
\end{equation}
results from a discrete group convolution over the grid, and the same relation gives $(Ku)_\ell = (\kappa * u)(z_\ell)$. The kernel $\kappa$ is thus the same for every pair of grid points and, in the barycenter problem, every marginal, so it is assembled once and reused throughout the iteration, as in the convolutional Wasserstein framework of \cite{solomon2015convolutional} and the iterative Bregman projections of \cite{benamou2015iterative}.

The convolutional structure identified above still leaves the kernel $\kappa$ to be numerically computed for all grid points. Evaluating $\kappa(z)=e^{-\mathring d_{\mathbb H,\lambda}(e,z)^p/\varepsilon}$ at a given offset $z$ requires $\mathring d_{\mathbb H,\lambda}(e,z)$, which is not available in fully closed form. The horizontal geodesic distance $\mathring d_{\mathbb H,\lambda}(e,z)$ is obtained by matching the enclosed horizontal area to the prescribed central displacement of $z$ in the characterization of Proposition~\ref{prop:geodesics}, resulting in a transcendental equation solved numerically for every offset. Assembling $\kappa$ therefore requires one such root-find per grid-point. 
Alternatively, we may use the Kor\'anyi distance $\mathring d_{K, \lambda}$ which admits the explicit expression \eqref{eq:koranyi_reduced} and \eqref{eq: koranyiintro}. It is not, however, a geodesic distance: unlike $\mathring d_{\mathbb H,\lambda}$, which is by construction an infimum over horizontal curve lengths, $\mathring d_{K,\lambda}$ does not arise from any length structure, and $(\mathring{\mathbb H},\mathring d_{K,\lambda})$ is not itself a geodesic space. Its relation to $\mathring d_{\mathbb H,\lambda}$, established in Proposition \ref{prop:cc-koranyi}, is
\begin{equation}
\mathring d_{K,\lambda}(e,z) \;\le\; \mathring d_{\mathbb H,\lambda}(e,z) \;\le\; \sqrt{\pi}\,\mathring d_{K,\lambda}(e,z), \qquad z\in\mathring{\mathbb H},
\end{equation}
with equality on the left on the horizontal plane $\sigma=0$ and equality up to the constant $\sqrt\pi$ on the central axis $t=\omega=0$ (Proposition~\ref{prop:cc-koranyi}). This bound transfers directly to the Gibbs kernel as follows:

\begin{lemma}\label{lem:kernel-sandwich}
Fix $p\ge1$ and $\varepsilon>0$, and let $K_{\mathrm{CC}}(x,y)=e^{-\mathring d_{\mathbb H,\lambda}(x,y)^p/\varepsilon}$ and $K_{K}(x,y)=e^{-\mathring d_{K,\lambda}(x,y)^p/\varepsilon}$ denote the Gibbs kernels built from the Carnot-Carathéodory and Korányi distances, respectively. Then, for all $x,y\in\mathring{\mathbb H}$,
\begin{equation}
K_K(x,y)^{(\sqrt\pi)^p} \;\le\; K_{\mathrm{CC}}(x,y) \;\le\; K_K(x,y).
\end{equation}
\end{lemma}

\begin{proof}
Since $\mathring d_{K,\lambda}\le \mathring d_{\mathbb H,\lambda}\le \sqrt\pi\, \mathring d_{K,\lambda}$, raising to the power $p\ge1$ gives $\mathring d_{K,\lambda}^p \le \mathring d_{\mathbb H,\lambda}^p \le (\sqrt\pi)^p \mathring d_{K,\lambda}^p$, and applying the decreasing map $t\mapsto e^{-t/\varepsilon}$ reverses the inequalities: $e^{-(\sqrt\pi)^p \mathring d_{K,\lambda }^p/\varepsilon} \le e^{-\mathring d_{\mathbb H,\lambda}^p/\varepsilon} \le e^{-\mathring d_{K,\lambda }^p/\varepsilon}$. The middle and right terms are $K_{\mathrm{CC}}$ and $K_K$, and the left term equals $K_K^{(\sqrt\pi)^p}$.
\end{proof}


\paragraph{Variable support.} Here the atoms $x_k,y_\ell$ are independent of one another and unstructured, so the offsets $x_k^{-1}\star y_\ell$ no longer range over a common finite difference set, so the kernel matrix $K_{k\ell}=e^{-\mathring d_{\mathbb H,\lambda}(x_k,y_\ell)^p/\varepsilon}$ cannot be used directly. Instead one may do a Lloyd-type alternation between the coupling and the positions themselves: given the current atoms, an optimal coupling $\pi^\star=(\pi^\star_{k\ell})$ solving the discrete transport (or barycenter) problem is computed, and the atoms are then displaced according to $\pi^\star$, along Carnot-Carathéodory geodesics, before the coupling is recomputed at the updated positions. This displacement step relies on the geodesics of $\mathring{\mathbb H}$ being available explicitly (Proposition~\ref{prop:geodesics}), so that each update is an explicit evaluation of a known curve rather than a numerical solve.

In the case of $n=2$ marginals barycenters are geodesics directly obtained from the knowledge of $\mathring d_{\mathbb H,\lambda}$ geodesics with no need for further iterations. Given $\mu,\nu$ as in \eqref{eq: disc_measures}, an optimal coupling $\pi^\star$ between them, and writing as $\gamma^t(x,y)$ a measurable selection of evaluation at time $t \in [0,1]$ constant-speed Carnot-Carathéodory geodesics between $x$ and $y$, the barycenter of $\mu,\nu$ with weights $(w_1,w_2)$, $w_1+w_2=1$, is obtained by McCann interpolation by pushing each atom $x_k$ a fraction $w_2$ of the way along its matched geodesic:
\begin{equation}
\mathrm{Bar}\big(\{\mu,\nu\},(w_1,w_2)\big) \;=\; (\gamma^{w_2})_\# \pi^\ast,
\end{equation}
which if $\mu,\nu$ have the same number of atoms $J=L$ in which case at least one optimal plan is given by a permutation $\sigma$, simplifies to
\begin{equation}
\mathrm{Bar}\big(\{\mu,\nu\},(w_1,w_2)\big) \;=\; \sum_{j=1}^J a_j\,\delta_{\gamma^{w_2}(x_j, y_{\sigma(j)})}.
\end{equation}
For barycenters with $n>2$ marginals such an explicit procedure of course no longer applies. Unlike in Euclidean space, in curved spaces such weighted Fr\'echet means cannot be computed just by iterations of geodesics, even if those are explicitly available between any two points. Indeed, Wasserstein space is known to be non-negatively curved in the sense of Alexandrov \cite[Sec.~7.3]{ambrosio2005gradient}.

\subsection{Optimality of right translation}\label{sec: right}
Right translations are not isometries of a left-invariant distance, but in some cases, they are still optimal Monge maps \cite{AMBROSIO2004}. 

\begin{proposition}\label{prop:right-translation}
Let $d\in\{\mathring d_{\Hb,\lambda},\mathring d_{K,\lambda}\}$, $p\ge1$, $\mu\in\mathcal P_p(\mathring\Hb)$, and let $h\in\mathring\Hb$ be such that $F_\lambda(h)=(s,\eta,0)$ is horizontal, equivalently $h=(\lambda s,\lambda^{-1}\eta,s\eta/2)$. Then $R_h$ is an optimal transport map from $\mu$ to $(R_h)_\#\mu$, and $W_p(\mu,(R_h)_\#\mu)=|(s,\eta)|$.
\end{proposition}

\begin{proof}
The projection $\mathrm{pr}$ is $1$-Lipschitz by Lemma \ref{lemma: projection}. Moreover $\mathrm{pr}(x\star h)=\mathrm{pr}(x)+v$, with $v\coloneq\mathrm{pr}(h)=(\lambda s,\eta/\lambda)$ and $|v|_{g_\lambda}=|(s,\eta)|$. Hence $d(0,h)\ge|(s,\eta)|$, with equality: for $\mathring d_{\Hb,\lambda}$, the curve $\xi\mapsto(\xi\lambda s,\xi\eta/\lambda,\xi^2s\eta/2)$, $\xi\in[0,1]$, is horizontal with velocity $\lambda sX_1+\lambda^{-1}\eta X_2$ of $\lambda$-norm $|(s,\eta)|$; for $\mathring d_{K,\lambda}$, $|F_\lambda(h)|_K=|(s,\eta)|$. By left invariance $d(x,x\star h)=d(0,h)$, so $(\mathrm{Id},R_h)_\#\mu$ has cost $|(s,\eta)|^p$, which is the cost of moving along straight lines. To prove the remaining lower bound, we project the measures to the plane to get
\begin{equation}
    W_p(\mu, (R_h)_\#\mu)\geq W_p^{(\R^2, g_\lambda)}(\operatorname{pr}_\#\mu, (\operatorname{pr}\circ R_h)_\#\mu).
\end{equation}
In particular, a simple computation shows that also the right-hand side has cost $|(s,\eta)|^p$.
\end{proof}

\section{Measure representations of complex valued signals}\label{sec: main}
The aim of this section is to associate with $f$ a probability measure on the reduced Heisenberg group $\rH$ which retains both the magnitude and the phase of its Gabor transform. We first define this lift and establish its equivariance under the Schr\"odinger representation (Sections \ref{subsec: well} and \ref{subsec: equi}). We then study reconstruction and Wasserstein stability (Section \ref{subsec: rec-stab}), and relate the horizontal geometry of the lifted phase surface to differential reassignment (Section \ref{sec: ass}).

\begin{assumption} \label{ass:M1} 
The window and the signal $\psi, f\in L^2(\R;\bC)$ satisfy $\|\psi\|_2=1$; $\tau\psi, \dot\psi\in L^2(\R,\bC)$; and $A_f^\psi\in L^1(\R^2)$.
\end{assumption}
These assumptions are natural in harmonic analysis: The $L^1$-integrability of the STFT amplitude guarantees $f,\psi\in M^{1,1}(\R)$ \cite[Chapter 11]{Grochenig2001}. Moreover, the regularity and the moment conditions on the window characterize the Sobolev-Shubin class 
\begin{equation}
    Q_1(\R):=\{\psi\in H^1(\R): (1+\tau^2)^{1/2} \psi\in L^2(\R)\}
\end{equation}
See \cite{Shubin01} or \cite[Lemma 8.25(e)]{Gro06}. Note that $\|\dot\psi\|_2>0$ for such windows.

We also set the scaling $\lambda>0$ of the metric in point \ref{it: scale} of Proposition \ref{prop:jet} and Lemma \ref{lemma: projection} to be window dependent in the following sense:
\begin{equation}\label{eq:lambda-choice}
    \lambda^2=\frac{\|\tau\psi\|_2}{\|\dot\psi\|_2}.
\end{equation}
For a centered window, this is the ratio between its time spread and its frequency spread. For example, in a real Gaussian window scaled as $\psi(\tau)\coloneq(\pi\lambda^2)^{-1/4} e^{-\tau^2/(2\lambda^2)}$ we recover $\lambda$ as above. This choice gives the lowest Lipschitz constant in the $L^2$ reconstruction (cf. Proposition \ref{eq:Rlip2}).

\subsection{Well-posedness of the phase lift}\label{subsec: well}
Introducing a cutoff parameter $\epsilon\geq 0$, we split the Gabor transform $G_\psi f$ into amplitude $A^\psi_f$ and phase $\Theta^\psi_f$ defined 
on 
\begin{equation}\label{eq: domain-eps}
    \Omega_{f, \epsilon}^\psi\coloneq\big\{(t,\omega)\colon A_f^\psi>\epsilon\|f\|_2\|\psi\|_2\big\}, \qquad \mathcal O_{f,\epsilon}^\psi\coloneq\Omega_{f,\epsilon}^\psi\times\mathbb T,
    \qquad
    \mathbb T\coloneq\R/2\pi\mathbb Z,
\end{equation}
with $\Theta^\psi_f$ as the argument of 
\begin{equation}
    u_f^\psi(t,\omega)\coloneq\frac{G_\psi f(t,\omega)}{|G_\psi f(t,\omega)|}\in S^1.
\end{equation}
We therefore define the phase graph directly as a map into the reduced
Heisenberg group,
\begin{equation}
\label{eq:phase-lift-map}
    \Gamma_f^\psi\colon\Omega_{f,0}^\psi\longrightarrow\rH,\qquad
    \Gamma_f^\psi(t,\omega)
    \coloneq(t,\omega,u_f^\psi(t,\omega)).
\end{equation}
Here we identify the phase variable of $\rH$ with $S^1\simeq
\R/2\pi\Z$. Whenever a local real-valued branch $\Theta_f^\psi$ is chosen,
$u_f=e^{i\Theta_f^\psi}$ and we write locally
\begin{equation}
    \Gamma_f^\psi(t,\omega)=
    (t,\omega,\Theta_f^\psi(t,\omega)).
\end{equation}
All differential identities involving $\Theta_f$ below are understood in this
local sense. They are independent of the chosen branch.

The restriction to positive $A^\psi_f$ included in $\Omega_{f,0}^\psi$ is unavoidable because the phase is undefined at the zeros of $G_\psi f$. More importantly, the phase becomes intrinsically ill-conditioned in a neighbourhood of those zeros, a phenomenon that is due to the pole-like behaviour of its derivative \cite{Balasz16}. 
Its gradient reads
\begin{equation} \label{eq:phase-gradient} 
\nabla\Theta_f^\psi = \operatorname{Im} \left( \frac{\nabla G_\psi f}{G_\psi f} \right). 
\end{equation}
Thus the phase map is uniformly Lipschitz on sets bounded away from zero, with a Lipschitz constant proportional to the inverse amplitude, but no uniform stability estimate can hold as the amplitude tends to zero.
Although the value of the phase at an exact zero is irrelevant for an $A_f^\psi\,dt\,d\omega$-weighted representation, points of very small magnitude can therefore carry a highly unstable phase. 
We suppress this poorly conditioned region by introducing an amplitude-dependent cutoff by increasing $\epsilon$.
\begin{definition}\label{def: meas}
For $\epsilon\in[0,1/2)$, let $\chi_\epsilon\colon [0,\infty)\to[0,1]$ be either the soft or hard cutoff 
\begin{equation}
  \chi_\epsilon^s(a)\coloneq\min\Bigl\{1,\ \max\bigl\{0,\ \tfrac{a-\epsilon}{\epsilon}\bigr\}\Bigr\},\qquad \chi_\epsilon^h(a)\coloneq\mathrm 1_{(\epsilon, +\infty)}(a),
\end{equation}
where $\chi_0(a)=1$ in either case.
For $f\in L^2(\R,\bC)\setminus\{0\}$, define the reliability factor 
\begin{equation}
    r_{f}^\psi(z)\coloneq\chi_\epsilon\!\Bigl(\tfrac{A_f^\psi(z)}{\|f\|_2\|\psi\|_2}\Bigr),
\end{equation}
and 
\begin{equation}
  w_f^\psi(z)\coloneq A_f^\psi(z)\,r_f^\psi(z),
  \qquad
  E_f^\psi\coloneq\int_{\R^2}w_f^\psi\,dz,\qquad \rho_f^\psi\coloneq\frac{w_f^\psi}{E_f^\psi},
\end{equation}
Unless explicitly stated otherwise, $\chi_\epsilon$ denotes either choice.
    We are assuming that $\epsilon$ is small enough so that $E_f>0$. In both Definition \ref{def: meas} and Equation \eqref{eq: domain-eps}, the norm $\|\psi\|_2$ appears for generality. In the following, we assume $\|\psi\|_2=1$ as mentioned in Assumption \ref{ass:M1}.
\begin{definition}
    For $f \in L^2(\R,\C)$ we define the lifted probability measure in $\calP(\mathring\Hb)$ as
    \begin{align}
    \mu_f^\psi\coloneq(\Gamma_f^\psi)_\#\bigl(\rho_f^\psi\,dz\bigr)\ \in\ \calP(\rH).
\end{align}
\end{definition}
The measure $\mu_f^\psi$ encodes the amount of rescaled energy in the amplitude $A_f^\psi=|G_\psi f|$ of $f$ with a threshold $\chi_\epsilon$, and `positioned' on the phase graph, which is the image of $\Gamma_f^\psi$. We will sometimes drop the superscript $\psi$ when there is no confusion in the chosen window $\psi$.
\end{definition}
Whenever $\epsilon>0$, the support of the measure is compact, and so $\mu_f^\psi\in \calP_p(\rH)$ for all $p\geq 1$. In the case $\epsilon=0$, this is not a guarantee, and the moment condition relates to the time-frequency decay in the standard modulation space $M^{1,1}_m(\R)$, see Remark \ref{rmk: moment mod}.
\begin{remark}
    The soft cutoff has the advantage that $a\mapsto a\chi_\epsilon^s(a/m)$ is Lipschitz, whereas the hard cutoff
introduces an artificial discontinuity at the threshold. For either cutoff,
for every $m>0$ and $a\in[0,m]$,
\begin{equation}
\label{eq:chidefect}
    0
    \leq
    a-a\chi_\epsilon(a/m)
    \leq
    \epsilon m\,
    \mathbf 1_{\{a\leq2\epsilon m\}}.
\end{equation}
Thus the cutoff affects only low-amplitude coefficients and discards at most a fraction $\epsilon m$ of amplitudes pointwise. 
\end{remark}

\subsection{Equivariance of the lift}\label{subsec: equi}
The essential structural property of $\mu_f$ is that the Schr\"odinger action on signals or windows becomes left or right translation on the reduced Heisenberg group.
\begin{lemma}\label{lemma: equi}
    Let $\Phi\in \{L^{-1}_h, R_h\}$ for some $h\in \rH$, and let $f, f', \psi, \psi'\in L^2(\R, \bC)\setminus \{0\}$ satisfy $\|f\|_2=\|f'\|_2$, $\|\psi\|_2=\|\psi'\|_2$ and $W_{\psi'} f'=W_\psi f\circ \Phi$. Then, for either choice of cutoff $\chi_\epsilon\in\{\chi_\epsilon^s,\chi_\epsilon^h\}$,
    \begin{equation}
        \mu_{f'}^{\psi'}=(\Phi^{-1})_\# \mu_f^\psi.
        \end{equation}
        \begin{proof}
            For $h=(s,\eta, \theta)$, write $\Phi(t, \omega, \phi)=(\tau(t, \omega), \phi+c(t,\omega))$, where $\tau, c$ also depend on $h$: this is clearly possible for $\Phi\in \{L^{-1}_h, R_h\}$ from the definition \eqref{eq:jet-group} of the group operation. Evaluating the hypothesis at $\phi=0$ gives $G_{\psi'}f'(t,\omega)=e^{-ic(t,\omega)}G_\psi f(\tau(t,\omega))$ for all $(t,\omega)$. Taking absolute values, $A^{\psi'}_{f'}=A^\psi_f\circ\tau$. Since $\tau$ is a translation and $\|f'\|_2=\|f\|_2$, this gives, for either cutoff, $r_{f'}^{\psi'}=r_f^\psi\circ\tau$, hence $w_{f'}^{\psi'}=w_f^\psi\circ\tau$ pointwise and, by translation-invariance of $L^2$, $E_{f'}^{\psi'}=E_f^\psi$ and $\rho_{f'}^{\psi'}=\rho_f^\psi\circ\tau$. 
            Taking arguments, $u_{f'}^{\psi'}=e^{-ic}\,(u_f^\psi\circ\tau)$ on $\Omega^{\psi'}_{f',0}=\tau^{-1}(\Omega^\psi_{f,0})$, so
        \begin{equation}
            \Phi\bigl(\Gamma_{f'}^{\psi'}(t,\omega)\bigr) =\bigl(\tau(t,\omega),\,e^{ic(t,\omega)}u_{f'}^{\psi'}(t,\omega)\bigr)=\bigl(\tau(t,\omega),\,u_f^\psi(\tau(t,\omega))\bigr)=\Gamma_f^\psi(\tau(t,\omega)).
        \end{equation}
        For $\zeta\in C_0(\rH)$, substituting $(t,\omega)=\tau^{-1}(t',\omega')$ (which has unit Jacobian), we get
        \begin{equation}
\int_{\rH}\zeta\,d\mu_{f'}^{\psi'}=\int_{\Omega^{\psi'}_{f',0}}\zeta(\Gamma_{f'}^{\psi'})\;\rho_{f'}^{\psi'}\dd t\,\dd\omega=\int_{\Omega^\psi_{f,0}}\zeta\bigl(\Phi^{-1}\Gamma_f^\psi(t',\omega')\bigr)\rho_f^\psi(t',\omega')\dd t'\dd\omega' =\int_{\rH}\zeta\circ\Phi^{-1}\,d\mu_f^\psi.  
        \end{equation}
        Hence the measures agree by duality.
        \end{proof}
\end{lemma}

\begin{theorem}\label{thm:lift-equivariance}
    Let $f,\psi\in L^2(\R,\bC)\setminus \{0\}$ be a signal and a window, respectively. It holds
    \begin{enumerate}
        \item\label{it: sc} \textbf{positive scaling:} for $c>0$, $\mu_{cf}=\mu_f$.
        \item\label{it:ph} \textbf{phase rotation:} for $\theta\in \R/2\pi\Z$,  $\mu_{P_\theta f}=(L_{(0,0,\theta)})_\# \mu_f$;
        \item\label{it: freq} \textbf{frequency modulation.} for $\eta\in \R$, $\mu_{M_\eta f} = (L_{(0,\eta,0)})_\# \mu_f$.
        \item\label{it: time} \textbf{time translation.} for $s\in \R$, $\mu_{T_s f} = (L_{(s,0,0)})_\# \mu_f$.
        \item \label{it: time-freq}\textbf{Schr\"odinger action on the signal.} for $h\in \rH$, $\mu_{U_h f} = \mu^\psi_{U_h f} = (L_{h})_\# \mu_f^\psi = (L_{h})_\# \mu_f$.
        \item \label{it: window}\textbf{Schr\"odinger action on the window.} for $h\in \rH$, $\mu_f^\psi= (R_h)_\#\mu_{f}^{U_h \psi}$.
    \end{enumerate}
    \begin{proof}
        For \ref{it: sc}, note that the positive scaling only affects the amplitude $A_{cf}=cA_f$ which is then rescaled by $E_{cf}=cE_f$ in the measure. Moreover, $r_{cf}=r_f$, so the cutoff is untouched. For all the other points, apply Lemma \ref{lemma: equi}.
    \end{proof}
\end{theorem}

\begin{corollary} \label{cor:wasserstein-invariance}
Let $d$ be a left-invariant metric on $\rH$, let $p\geq1$, and assume
$\mu_f,\mu_g\in\mathcal P_p(\rH)$. Then for every $h\in\rH$,
\begin{equation}
\label{eq:W-invariance}
    W_{p}\bigl(\mu_{U_hf},\mu_{U_hg}\bigr)=W_{p}(\mu_f,\mu_g).
\end{equation}
\end{corollary}

\begin{proof}
By Theorem~\ref{thm:lift-equivariance},
\begin{equation}
    \mu_{U_hf}=(L_h)_\#\mu_f,
    \qquad
    \mu_{U_hg}=(L_h)_\#\mu_g.
\end{equation}
The claim follows from the left invariance of $d$ by Proposition \ref{prop:Wp-left-invariant}.
\end{proof}

\begin{corollary}\label{cor: monge}
    Let $d$ be either Carnot-Carath\'eodory or the Kor\'anyi distance on $\rH$, let $p=2$, and $h=(\lambda s,\eta/\lambda, s\eta/2)$. Then $R_h$ is an optimal Monge map between $\mu_f^{U_h\psi}$ and $\mu_f^\psi$ for $d^2$.
    \begin{proof}
        It follows from Point \ref{it: window} of Theorem \ref{thm:lift-equivariance}, together with the previous discussion of optimality of right translation for left-invariant metrics in Proposition \ref{prop:right-translation}.
    \end{proof}
\end{corollary}

\begin{remark}\label{rmk: radar}
    A possible use of the equivariance result lies in radar signal processing. With the transmitted pulse $f$ as window, $G_f r(s,\eta)=\langle r,T_sM_\eta f\rangle$ is the cross-ambiguity function of a received signal $r$ and $f$.
    Its amplitude $A^f_r$ is the matched-filter output, whose maximizer estimates delay and Doppler shift \cite{Woodward1953}, \cite[Section 1.4]{Folland89}. The transform $W_f r$ extends it to $\mathring\Hb$, as in the group-theoretic treatment of ambiguity functions in \cite{Aus85}. 
    By Equation \eqref{eq: full-gabor}, 
    \begin{equation}
        \operatorname{Re}W_f r(s,\eta,\theta)=A^f_r(s,\eta)\cos\big(\Theta^f_r(s,\eta)-\theta\big),
    \end{equation}
    so $\max_{\mathring\Hb}\operatorname{Re}W_f r=\max_{\R^2}A^f_r$, and wherever $A^f_r>0$ the maximizing phase is unique, $\theta=\Theta^f_r(s,\eta)$. Phase-coherent matching therefore has the value of matched filtering, and its maximizers lie on the phase graph $\Gamma^f_r$, which carries $\mu^f_r$.
    For an echo $r=c\,U_{h^*}f$ with $c>0$, Lemma \ref{lem:Wcov} gives $W_f r=c\,W_f f\circ L_{h^*}^{-1}$. 
    Its real part is maximal exactly at $h^*$. Theorem \ref{thm:lift-equivariance} gives $\mu^f_r=(L_{h^*})_\#\mu^f_f$: in the lift with the pulse as window, an echo is the left translate of the lifted auto-ambiguity function by its delay, Doppler shift and phase.
\end{remark}

\subsection{Reconstruction and Wasserstein inverse stability}\label{subsec: rec-stab}
The lift $\mu_f$ is a positive measure and therefore no longer carries the
complex STFT coefficients directly. Nevertheless, it retains enough
information to reconstruct the part of the Gabor expansion which survives
the amplitude cutoff.
\begin{proposition}\label{prop: rec}
    Let $\mathcal{R}_\psi:\calP(\rH)\to L^\infty(\R)$ be the linear operator
\begin{equation}\label{eq:Rdef}
  \mathcal{R}_\psi\mu(\tau)\coloneq\frac{1}{2\pi}\int_{\rH} U_h\psi \,d\mu(h).
\end{equation}
Under Assumption \ref{ass:M1}, with $\|f\|_2=1$
and $\tilde f\coloneq E_f\,\mathcal{R}_\psi\mu_f$,
\begin{equation}\label{eq:recon}
  \bigl\|f-\tilde f\bigr\|_{\infty}\leq \frac{\|\psi\|_\infty}{2\pi}\int_{\{A_f\leq 2\epsilon\}} A_f(z)  \dd z.
\end{equation}

\begin{proof}
    Since $E_f\rho_f=w_f$ and $w_fe^{i\Theta_f}=r_fA_fe^{i\Theta_f}
=\chi_\epsilon(A_f)\,G_\psi f$,
\begin{equation}
  \tilde f(\tau)=\frac{1}{2\pi}\int_{\R^2}\chi_\epsilon (A_f(t,\omega))\,G_\psi f(t,\omega)\,
    e^{i\omega(\tau-t)}\psi(\tau-t)\dd\omega\dd t .
\end{equation}
Subtracting from the Gabor inversion formula \eqref{eq:gabor-inversion}
\begin{equation}
    f(\tau)=\frac{1}{2\pi} \int_{\R^2} G_\psi f(t,\omega)\,e^{i\omega(\tau-t)}\psi(\tau-t)\,\dd\omega\dd t, 
\end{equation}
we get,
\begin{equation}
  \bigl|f(\tau)-\tilde f(\tau)\bigr|
  \le\frac{\|\psi\|_\infty}{2\pi}\int_{\R^2}\bigl(A_f-A_f\chi_\epsilon(A_f)\bigr)\dd z .
\end{equation}
By Lemma \ref{lemma: AC0} and the assumption on unitary $2$-norm, $0\leq A_f^\psi\leq \|f\|_2\|\psi\|_2=1$. Therefore, using the pointwise estimate in \eqref{eq:chidefect} with $m=1$, we get the claim.
\end{proof}
\end{proposition}
Since $A_f\in L^1$, the integral is finite and tends to $0$ as $\epsilon\to 0$ by dominated convergence.

\begin{corollary}\label{cor: rec1}
    For $\epsilon=0$, the pair $(\mu_f^\psi, \psi)$ identifies $f$ uniquely up to positive scaling.
    \begin{proof}
        Assume $f,g\neq 0$ such that $\mu_f=\mu_g$. Set $\hat{f}, \hat{g}$ to be the $L^2$ normalized signals. By Theorem \ref{thm:lift-equivariance}, $\mu_{\hat{f}}=\mu_{\hat{g}}$. In the absence of cutoff, i.e. $\epsilon=0$, Proposition \ref{prop: rec} yields
        $\hat{f}=E_{\hat{f}}\calR_\psi\mu_{\hat{f}}$. Likewise $g=E_{\hat{g}}\calR_\psi\mu_{\hat{g}}$, which proves that $\hat{g}$ and $\hat{f}$ are multiple of each other, $\hat{g}=\left(E_{\hat{g}}/E_{\hat{f}}\right)\hat{f}$.
    \end{proof}
\end{corollary}

\begin{corollary}
    Let $\epsilon=0$, fix a window $\psi$, and let 
    \begin{equation}
        \calA_{\psi,p}\coloneq\{f\in L^2(\R, \bC)\setminus \{0\}\colon \mu_f^\psi\in \calP_p(\rH)\}.
    \end{equation}
    Let $f\sim_+ g$ if $g=cf$ for some $c>0$. Then, 
    \begin{equation}
        d_{\psi, p} ([f]_+, [g]_+)\coloneq W_p(\mu_f^\psi, \mu_g^\psi)
    \end{equation}
    defines a metric on $\calA_{\psi, p}/\sim_+$. Moreover, for every $h\in \rH$, 
    \begin{equation}
        d_{\psi, p}([U_h f]_+, [U_hg]_+)=d_{\psi, p} ([f]_+, [g]_+).
    \end{equation}
    \begin{proof}
        Positive-scale invariance follows from Theorem \ref{thm:lift-equivariance}, while separation follows from Corollary \ref{cor: rec1}. The remaining metric properties are inherited from $W_p$ and Schr\"odinger invariance follows from Corollary \ref{cor:wasserstein-invariance}.
    \end{proof}
\end{corollary}

\begin{remark}\label{rmk: moment mod}
    A reasonable question to ask at this point is whether the $p$-moment condition on the measure $\mu_f^\psi$ relates to the regularity of $f$ in the modulation spaces $M^{q,r}_m(\R)$ \cite[Chapter 11]{Grochenig2001}.  
    Let $\psi\in \calS(\R)\setminus \{0\}$ be of Schwartz class, and let $\epsilon=0$. 
    For $p\geq 1$, define the polynomial weight 
    \begin{equation}
        m_p(z):=\left(1+|z|_{g_\lambda}\right)^p. 
    \end{equation}
    The space $\calA_{\psi,p}$ coincides with $M^{1,1}_{m_p}(\R)\setminus \{0\}$. Indeed, Lemma~\ref{lemma: projection} and the compactness of $\mathbb T$ yield
\begin{equation}
    |z|_{g_\lambda}
    \leq \mathring d_{\mathbb H,\lambda}(0,(z,\phi))
    \leq |z|_{g_\lambda}+2\pi,
\end{equation}
so that the finite $p$-moment condition is equivalent to
\begin{equation}
    \int_{\R^2}m_p(z)|G_\psi f(z)|\,dz<\infty.
\end{equation}
In contrast, inclusions of the form $M^{q,q}(\R)\subset\calA_{\psi,p}$, for $q\geq1$, are generally false. Thus the moment condition affects the polynomial weight rather than the integrability exponents of the modulation space.
In particular, increasing $p$ strengthens the required time-frequency decay while leaving the integrability exponents fixed at $1$ by definition of $\mu_f^\psi$.
\end{remark}


We can push the reconstruction theory even further with a cross-window result, useful in view of the right translation of the measure as Schr\"odinger action on the window, see point \ref{it: window} of Theorem \ref{thm:lift-equivariance}.
\begin{lemma}\label{lemma: rec2}
    For $\epsilon=0$ and $\psi'\in L^2(\R)\setminus \{0\}$, the linear synthesis operator $\calR_{\psi'}$ reconstructs $f$ uniquely up to scaling from $\mu_f^\psi$ if and only if $\langle \psi, \psi'\rangle\neq 0$.
    \begin{proof}
        Note that by definition of Schr\"odinger representation and phase lift, 
        \begin{equation}
            U_{\Gamma_f^\psi(t,\omega)} \psi'=e^{i\Theta_f^\psi(t,\omega)} U_{(t,\omega, 0)}\psi'.
        \end{equation}
        Therefore, for $C_{\psi'}=2\pi\|\psi'\|^2_2$,
        \begin{align}
            \calR_{\psi'}\mu_f^\psi=\frac{1}{C_{\psi'}}\int_{\rH} U_h\psi'\dd\mu_f^\psi= \frac{1}{E_f^\psi C_{\psi'}}\int_{\R^2} G_\psi f(t,\omega) U_{(t,\omega, 0)}\psi'\dd t\dd \omega.
        \end{align}
        For any $g\in L^2(\R, \bC)$, apply Moyal's identity,  
        \begin{align}
            \langle \calR_{\psi'}\mu_f^\psi, g\rangle &= \frac{1}{E_f^\psi C_{\psi'}}\int_{\R^2} G_\psi f(t,\omega) \langle U_{(t,\omega, 0)}\psi',g \rangle \dd t\dd \omega\\
            &=\frac{1}{E_f^\psi C_{\psi'}}\int_{\R^2} G_\psi f(t,\omega) \overline{G_{\psi'}g(t,\omega)} \dd t\dd \omega\\
            &=\frac{2\pi}{E_f^\psi C_{\psi'}} \langle \psi', \psi\rangle\langle f, g\rangle.
        \end{align}
        Therefore, 
        \begin{equation}
            E_f^\psi \calR_{\psi'}\mu_f^\psi= \frac{\langle \psi', \psi\rangle}{\|\psi'\|_2^2} f.
        \end{equation}
    \end{proof}
\end{lemma}

\begin{proposition}\label{prop:Rlip2}
Let $\psi$ be real-valued. The operator $\mathcal{R}_\psi$ of \eqref{eq:Rdef} satisfies, for all $\mu,\nu\in\calP_1(\rH)$,
\begin{equation}\label{eq:Rlip2}
  \bigl\|\mathcal{R}_\psi\mu-\mathcal{R}_\psi\nu\bigr\|_{2}\leq \,L^{(2)}_\psi\;W_1(\mu,\nu),
  \qquad L^{(2)}_\psi=\frac{\sqrt{\|\dot\psi\|_2\|\tau\psi\|_2}}{2\pi}.
\end{equation}
Moreover, if $\psi\in W^{1,\infty}(\R)$ and $\tau\psi\in L^\infty(\R)$, it holds
\begin{equation}\label{eq:Rlip}
  \bigl\|\mathcal{R}_\psi\mu-\mathcal{R}_\psi\nu\bigr\|_{\infty}\leq \,L^{(\infty)}_\psi\;W_1(\mu,\nu),
  \qquad L^{(\infty)}_\psi=\frac{\bigl(\lambda^2\|\dot \psi\|_\infty^2+\lambda^{-2}\|\tau\psi\|_\infty^2\bigr)^{1/2}}{2\pi}.
\end{equation}
\begin{proof}
Let $q\in\{2,\infty\}$ and $L=L^{(q)}_\psi$. We first show
\begin{equation}\label{eq:orbit-lip}
  \|U_x\psi-U_{y}\psi\|_q\le 2\pi L\,\mathring d_{\Hb,\lambda}(x,y),\qquad x,y\in\rH.
\end{equation}
Let $\gamma\colon[0,1]\to\rH$ be a horizontal curve from $x$ to $y$, with $\dot\gamma=aX_1^\lambda+bX_2^\lambda$. By Lemma~\ref{lemma:U-deriv}, for a.e.\ $\xi$,
\begin{equation}
    \tfrac{d}{d\xi}U_{\gamma(\xi)}\psi=U_{\gamma(\xi)}g_\xi,\qquad g_\xi\coloneq -a(\xi)\lambda\dot\psi+ib(\xi)\lambda^{-1}\tau\psi .
\end{equation}
This is an $L^2$-valued derivative, which covers the case $q=2$. For $q=\infty$, the same computation holds pointwise in $\tau$, since $\psi$ is Lipschitz.
Every $U_h$ is an isometry of $L^q$, so Minkowski's integral inequality gives
\begin{equation}
  \|U_{y}\psi-U_x\psi\|_q \le \int_0^1\|U_{\gamma(\xi)}g_\xi\|_q\dd\xi =\int_0^1\|g_\xi\|_q\dd\xi .
\end{equation}
For $q=2$ we can use that $\psi$ is real, so $\dot\psi$ and $i\tau\psi$ are pointwise orthogonal in $\bC$. By the choice \eqref{eq:lambda-choice} of $\lambda$,
\begin{equation}
  \|g_\xi\|_2^2=a^2\lambda^2\|\dot\psi\|_2^2+b^2\lambda^{-2}\|\tau\psi\|_2^2=(a^2+b^2)\,\|\dot\psi\|_2\|\tau\psi\|_2 .
\end{equation}
For $q=\infty$, the Cauchy-Schwarz inequality gives
\begin{equation}
    \|g_\xi\|_\infty\le|a|\lambda\|\dot\psi\|_\infty+|b|\lambda^{-1}\|\tau\psi\|_\infty
  \le(a^2+b^2)^{1/2}\bigl(\lambda^2\|\dot\psi\|_\infty^2+\lambda^{-2}\|\tau\psi\|_\infty^2\bigr)^{1/2}.
\end{equation}
In both cases $\|g_\xi\|_q\le 2\pi L\,|\dot\gamma(\xi)|_\lambda$. Minimizing the length of $\gamma$ yields~\eqref{eq:orbit-lip}.
Now let $\pi\in\Pi(\mu,\nu)$. By~\eqref{eq:Rdef}, Minkowski's integral inequality and \eqref{eq:orbit-lip},
\begin{equation}
  \|\mathcal R_\psi\mu-\mathcal R_\psi\nu\|_q
  =\frac1{2\pi}\Bigl\|\int_{\rH\times\rH}(U_x\psi-U_{y}\psi)\dd\pi(x,y)\Bigr\|_q
  \le L\int_{\rH\times\rH}\mathring d_{\Hb,\lambda}\,d\pi .
\end{equation}
Taking the infimum over $\pi$ concludes the proof.
\end{proof}
\end{proposition}
This Lipschitz continuity can be adapted to the Kor\'anyi metric by adjusting the Lipschitz constant. Moreover, it extends to $W_p$ immediately for $p\geq1$.


\subsection{Differential reassignment as a Wasserstein gradient flow}\label{sec: ass}
The interpretation of differential reassignment as a type of evolution equation in the Heisenberg group was extensively studied in \cite{DUITS2013}. Here, we propose a Wasserstein gradient flow viewpoint based on our measure construction.

\begin{lemma}
    Assume $G_\psi f \in C^1(\Omega_{f,0})$. Then, the phase surface $\Sigma_f\coloneq\Gamma_f(\Omega_{f,0})$ is $C^1$, with $T\Sigma_f=\ker(\mathrm d\Theta_f-\mathrm d\phi)$. The horizontal field
\begin{equation}\label{eq:Vf}
  V_f\coloneq(\hat t_f-t)\,X_1+(\hat\omega_f-\omega)\,X_2,
  \qquad
  \hat\omega_f\coloneq\partial_t\Theta_f,\quad \hat t_f\coloneq t-\partial_\omega\Theta_f,
\end{equation}
spans $T\Sigma_f\cap D$ wherever $V_f\neq0$; at the characteristic points
$\{V_f=0\}$ one has $T\Sigma_f= D$.
\begin{proof}
    The surface $\Sigma_f\coloneq\Gamma_f(\Omega_{f,0})$ can equivalently be characterized as
\begin{equation}
    \Sigma_f=\left\{\operatorname{Im}(W_\psi f)^{-1}(\{0\})\right\}\cap \{\operatorname{Re}(W_\psi f)>0\},
\end{equation}
By assumption, $W_\psi f$ is $C^1$, and on $\Sigma_f$
\begin{equation}
    \partial_\phi \operatorname{Im}(W_\psi f)=\operatorname{Im}(-i W_\psi f)=-\operatorname{Re}(W_\psi f)\neq 0.
\end{equation}
By the implicit function Theorem, we get that $\Sigma_f$ is $C^1$. Locally, after choosing a smooth phase branch $\Theta_f$,
\begin{equation}
\label{eq:W-ker}
    T \Sigma_f=\ker(d\Theta_f-d\phi).
\end{equation}
The other claims follow from the fact that $X_1$ and $X_2$ span $D$, and that 
        \begin{align}
            V_f&=\left\{\left[d\Theta_f(X_1)-d\phi(X_1)\right] X_2- \left[d\Theta_f(X_2)-d\phi(X_2)\right] X_1\right\}\\
            &=\{\left[-\partial_\omega \Theta_f\right] X_1+\left[\partial_t\Theta_f -\omega\right] X_2\}.\qedhere
        \end{align}
\end{proof}
\end{lemma}

The coefficients of $V_f$ do not depend on $\phi$, so $V_f$ is defined on
$\mathcal O_{f,0}$. Its planar part $v_f\coloneq(\hat t_f-t,\hat\omega_f-\omega)$ is the displacement of differential reassignment \cite{Cha97, daudet99}, and the integral curves of $V_f$ are horizontal curves preserving the Gabor phase \cite[Lemma 9]{DUITS2013}. Classical reassignment \cite{Kodera1976,Kodera1978,Auger1995} applies $z\mapsto z+v_f(z)$ once; synchrosqueezing \cite{daubechies2009, Thakur_2011} applies its frequency component only. 
For Gaussian windows, all three are gradient dynamics.
Let $\psi$ be a fixed real-valued Gaussian window with $\|\psi\|_2=1$. More precisely, $\psi=\psi_\lambda(\tau)=(\pi\lambda^2)^{-1/4} e^{-\tau^2/(2\lambda^2)}$, for $\lambda>0$ the window parameter determining the scaling in point \ref{it: scale} of Proposition \ref{prop:jet} and Lemma \ref{lemma: projection}.

\begin{proposition}\label{prop:gradient}
Let $\Psi_f\coloneq\log A_f\circ\operatorname{pr}$ on $\mathcal{O}_{f,0}$. Then
\begin{equation}\label{eq:grad}
  V_f=\nabla_{\rH, \lambda}\Psi_f,
  \qquad\text{equivalently}\qquad
  v_f=\nabla_{g_\lambda}\log A_f
     =\begin{pmatrix}\lambda^2\partial_t\log A_f\\ \lambda^{-2}\partial_\omega\log A_f\end{pmatrix},
\end{equation}
where $\nabla_{\rH, \lambda}$ and $\nabla_{g_{\lambda}}$ are the subRiemannian and the flat $(\R^2, g_\lambda)$ gradients, respectively.
In particular $|V_f|_\lambda\circ\Gamma_f=|v_f|_{g_\lambda}$, and the characteristic points
of $\Sigma_f$ are the lifts of the critical points of $A_f$.
\begin{proof}
Take $n=1$ and $c=0$ in the setting of \cite{DUITS2013}, and set $a\coloneq\sqrt{2\pi}\,\lambda$. The isomorphism
$(p,q,s)\mapsto(p,2\pi q,2\pi s+\pi pq)$ identifies their Schr\"odinger representation and
Gabor transform $\mathcal W_{\psi_a}f$ with $U$ and $W_{\psi_\lambda}f$.
It is the composition of $F_1^{-1}$ and $(s,\eta, \theta)\mapsto (s, 2\pi \eta, 2\pi\theta)$, so, using Proposition \ref{prop:equi}, it also maps their left-invariant fields $\mathcal A_1,\mathcal A_2$ to $X_1,2\pi X_2$, and their phase
$\Omega=\arg\mathcal W_{\psi_a}f$ to $\Theta_f-\phi$. The Cauchy-Riemann relations of \cite[(33), (36)]{DUITS2013} then read
\begin{equation}
    \partial_t\Theta_f-\omega=X_1(\Theta_f-\phi)=\lambda^{-2}\partial_\omega\log A_f,\qquad
  \partial_\omega\Theta_f=X_2(\Theta_f-\phi)=-\lambda^{2}\partial_t\log A_f .
\end{equation}
Substituting into $V_f=-\partial_\omega\Theta_f\,X_1+(\partial_t\Theta_f-\omega)\,X_2$, gives
\begin{equation}
    V_f=\lambda^2\partial_t\Psi_f\,X_1+\lambda^{-2}\partial_\omega\Psi_f\,X_2=X_1^\lambda(\Psi_f)X_1^\lambda +X_2^\lambda(\Psi_f)X_2^\lambda=\nabla_{\rH, \lambda}\Psi_f,
\end{equation} since
$\Psi_f$ does not depend on $\phi$. The remaining claims follow from
$|aX_1+bX_2|_\lambda^2=\lambda^{-2}a^2+\lambda^2b^2$ and $A_f>0$ on $\mathcal O_{f,\epsilon}$.
\end{proof}
\end{proposition}

Proposition~\ref{prop:gradient} is the subRiemannian form of the observation of \cite{DUITS2013} that, for Gaussian windows, the reassignment field derives from the log-spectrogram: the reassignment direction is the unique horizontal direction tangent to the phase surface which when projected to $\R^2$ results in the steepest-ascent direction of $\log A_f$. The $\lambda$ scaling in the subRiemannian metric makes the gradient structure particularly clean.

Extend $\Psi_f$ by $-\infty$ outside of $\mathcal{O}_{f,0}$, and define the potential energy
\begin{equation}
    \mathcal V_f(\mu)\coloneq -\int_{\rH}\Psi_f\,\mathrm d\mu,
  \qquad \mu\in\mathcal P_2(\rH),
\end{equation}
which measures the (negative) expected log-amplitude seen by $\mu$. 

\begin{proposition}\label{prop: non-conv}
    For all $f\in L^2(\R,\bC)\setminus\{0\}$, the energy $\mathcal{V}_f$ is not geodesically convex in $(\calP(\rH), W_2)$, where $W_2$ is induced by $\mathring d_{\Hb,\lambda}$. 
    More precisely, there exist $x,y\in \mathcal{O}_{f,0}$ joined by a unique subRiemannian geodesic $\gamma\colon [0,1]\to \rH$ such that 
    \begin{equation}
        \mathcal{V}_f(\mu_{1/2})- \frac{1}{2}\mathcal{V}_f(\mu_0)- \frac{1}{2}\mathcal{V}_f(\mu_1)>0,
    \end{equation}
    where $\mu_\xi\coloneq\delta_{\gamma(\xi)}$.
\end{proposition}
\begin{proof}
On Dirac masses, $\mathcal{V}_f(\delta_x)=-\log A_f(\operatorname{pr }x),$ so convexity along $\delta_{\gamma(\xi)}$ just asks that $-\log A_f$ be convex along the time-frequency projection of every Heisenberg geodesic. Those projections are not straight lines. In the rescaled plane $z=\Lambda(t,\omega)=(t/\lambda, \lambda\omega)$, which consists of the first two coordinates of $F_\lambda$, they are circular arcs and the central difference between the endpoints in $F_\lambda$, which equals the enclosed area, sets how much the arc bends and to which side. In the original $(t,\omega)$ variables, they are arcs of ellipses with equal enclosed area because $\Lambda$ has determinant $1$.
Since $\mathcal{V}_f$ ignores phase, we can choose any phase difference for free, and in particular force the geodesic to project onto half-ellipses bulging in a direction of our choice.
Now set $u\coloneq\log A_f\circ \Lambda^{-1}$ and take $z\in \Lambda(\Omega_{f,0})$ with $\nabla u(z)\neq 0$ and place two endpoints $z_0,z_1$ symmetrically at Euclidean distance $\ell/2$ from $z$ on the line orthogonal to the gradient, and make the half-circle bulge in the direction $-\nabla u(z)$. We can make $\ell$ small enough so that the phase difference remains between $0$ and $\pi$, i.e. $\ell<2\sqrt{2}$, and small enough so that the closed disk of radius $\ell/2$ lies in $\Lambda(\Omega_{f,0})$. 
Take the Taylor expansion of $u$ at the point $z$:
\begin{equation}
    u(z+v)=u(z)+\nabla u(z)\cdot v + R(v),
\end{equation}
for some first order residual $R$ satisfying
\begin{equation}
    o_R(\rho)\coloneq\sup_{0<|v|\leq \rho}\frac{|R(v)|}{|v|}\xlongrightarrow{\rho\to 0} 0.
\end{equation}
Let $n\coloneq\nabla u(z)/|\nabla u(z)|$ and $e$ be the unit vector orthogonal to $n$. Apply the expansion with $v=\pm \frac{\ell}{2}e$ and $v=-\frac{\ell}{2}n$: first note that for all of the choices $|R(v)|\leq \frac{\ell}{2} o_R (\ell/2)$, then
\begin{equation}
    u\left(z-\frac{\ell}{2}n\right) -\frac{1}{2}(u(z_0)+u(z_1))\leq -\frac{\ell}{2}|\nabla u(z)| + \ell o_R\left(\frac{\ell}{2}\right).
\end{equation}
Again, by choosing $\ell$ small enough, we can ensure that $o_R(\ell/2)<\frac{1}{2}|\nabla u(z)|$, which gives the desired inequality.  Being $z-\ell/2 n$ the mid point of the circular arc, we can reapply $\Lambda^{-1}$, and get the desired unique geodesic.
\end{proof}

Formally, the first variation of $\mathcal{V}_f$ is $-\Psi_f$. Hence, the formal $W_2$-gradient flow velocity field associated with $\mathcal{V}_f$ is 
\begin{equation}
    -\nabla_{\rH,\lambda}\frac{\delta\mathcal{V}_f}{\delta\mu}=\nabla_{\rH,\lambda}\Psi_f=V_f.
\end{equation}
Thus differential reassignment has the formal continuity-equation representation
\begin{equation}
    \partial_s \mu_s +\operatorname{div}_{\rH,\lambda}(\mu_s V_f)=0.
\end{equation}


\section{Numerical experiments}\label{sec: num}
Throughout the experiments we work on the reduced Heisenberg group $\mathring{\mathbb H}$, using a uniform discretization of the time-frequency domain and periodic discretization of the phase variable. The geometric scale $\lambda$ is chosen consistently with the Gaussian analysis window $\psi_\lambda$; equivalently, the grid is scaled so that $\Delta t/\lambda$ and $\lambda\Delta\omega$ have comparable resolution. Unless stated otherwise, Wasserstein computations use the scaled Heisenberg distance $\mathring d_{\mathbb H,\lambda}$. For the signal interpolation experiments we use barycentric parameters
\begin{equation}
    \vartheta\in\{0,0.25,0.5,0.75,1\}.
\end{equation}
The phase and the corresponding time-frequency intensity are displayed separately in order to visualize both the evolution of the lifted phase structure and its projected signal representation.

For the interpolation $\mu_\vartheta$, we also estimate how much $\mu_\vartheta$ represents a graph measure through a graph index $\mathcal G$. For $\mu\in\mathcal P(\mathring H)$ with time–frequency marginal $\rho$ and disintegration $\mu=\int\mu_z d\rho(z)$, let 
\begin{equation}
    \mathcal G(\mu)\coloneq\int_{\mathbb R^2}\left|\int_{\mathbb T}e^{i\theta} d\mu_z(\theta)\right| d\rho(z)\in[0,1].
\end{equation} 
Since $\mu_z$ is a probability measure on $\mathbb T$, its mean $m(z)=\int_{\mathbb T} e^{i\theta}\dd \mu_z(\theta)$ satisfies $|m(z)|\leq 1$, with equality if and only if $\mu_z$ is a Dirac mass. Hence, $\mathcal G(\mu)=1$ if and only if $\mu$ is carried by the graph of a phase function. $\mathcal G$ is left-invariant, since left translations shift the phase by a constant on each fibre. 
For the discrete $\mu_\vartheta$ we coarse-grain on the cells $(Q_j)_j$ of the STFT grid, $\mathcal G_{Q}(\mu)=\sum_{j}\big|\int_{Q_j\times\mathbb T}e^{i\tilde\theta} d\mu\big|$, where $\tilde\theta=\theta-\bar\omega(t-\bar t)-\frac12(t-\bar t)(\omega-\bar\omega)$ is the central coordinate of $c_\mu^{-1}\star x$ and $c_\mu=(\bar t,\bar\omega,0)$ is the centroid of $\rho$. 
Since the endpoints $\mu_0,\mu_1$ are lifts, we have $\mathcal G_{Q}(\mu_0)=\mathcal G_{Q}(\mu_1)=1$. 
In the Lagrangian setting of Section \ref{subsec: McCann}, the interpolators consist of off-the-grid particles whose projected locations are generically distinct, so the exact disintegration of the atomic measure would yield a Dirac conditional phase at each location and therefore $\mathcal G(\mu)=1$ trivially.
Instead, for $\mu_\vartheta=\sum_{k}a_k \delta_{x_k}$ with $x_k=(z_k,\theta_k)$,
we estimate the conditional $m(z)=\mathbb E[e^{i\theta}\mid z]$ by the Nadaraya–Watson estimator
\begin{equation}
    \widehat m_b(z):= \frac{\sum_k a_k K_b(z-z_k) e^{i\widetilde \theta_k(z)}}{\sum_{k} a_k K_b(z-z_k)}
\end{equation}
using the same locally centered phase $\widetilde \theta_k(z)$ as in the cell-based definition above, and the positive kernel on the time-frequency plane with bandwidth $b$ (consistent with the $Q_j$)
\begin{equation}
    K_b(z,z_k):=\exp\left[-\frac{|z-z_k|_{g_\lambda}^2}{2b^2}\right].
\end{equation}
We then estimate the graph index from $|\widehat m_b|$ at the STFT resolution $b$.
In the Eulerian setting of Section \ref{subsec: bary}, the measures are discretized on a group-consistent grid, on which we also estimate the graph index $\calG$.


\begin{figure}[!ht]
    \centering
    \includegraphics[width=1\linewidth]{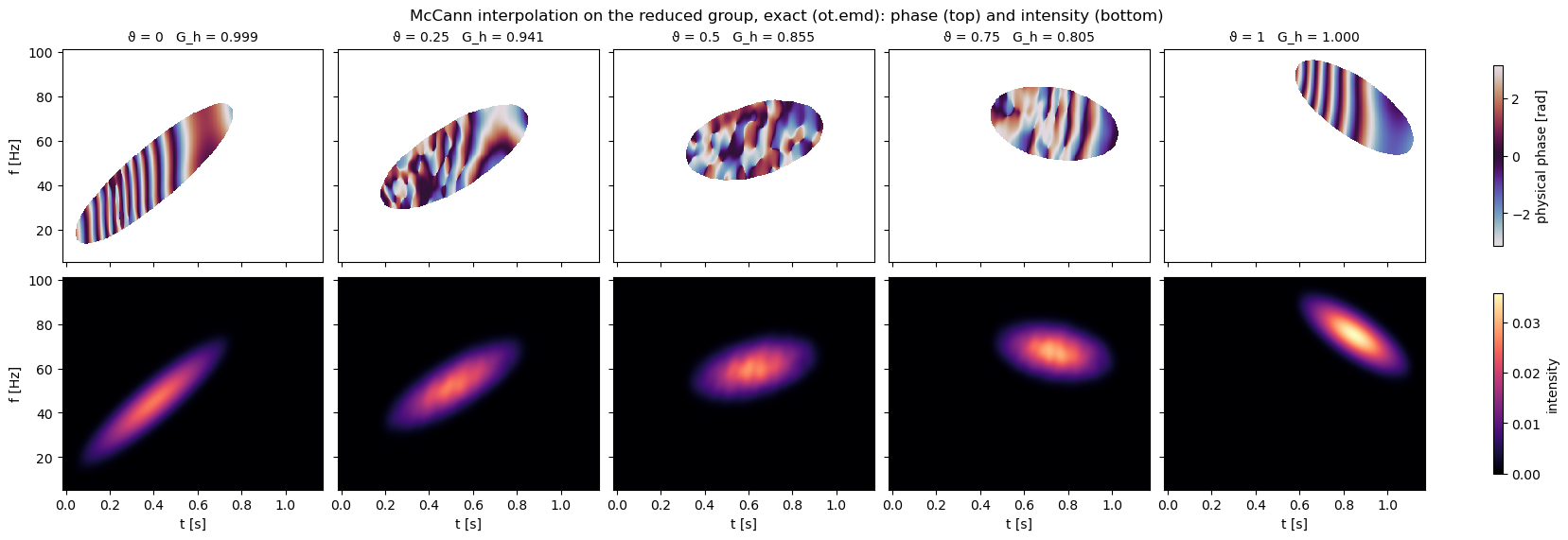}
    \includegraphics[width=1\linewidth]{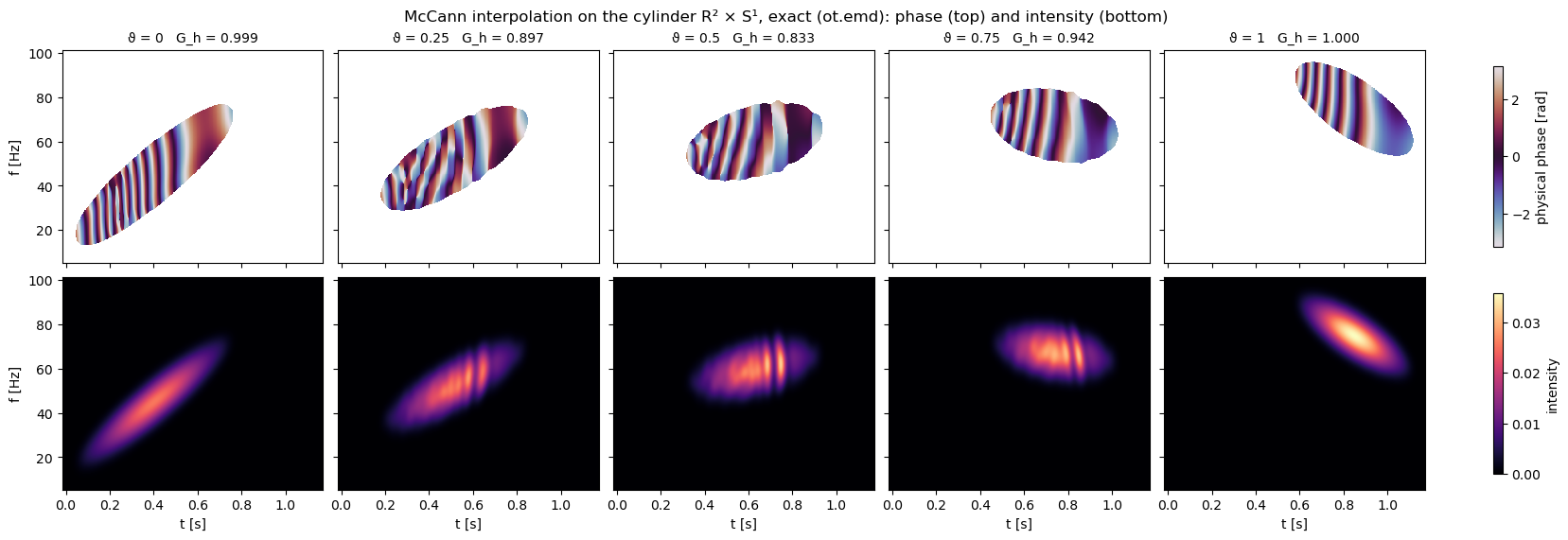}
    \includegraphics[width=1\linewidth]{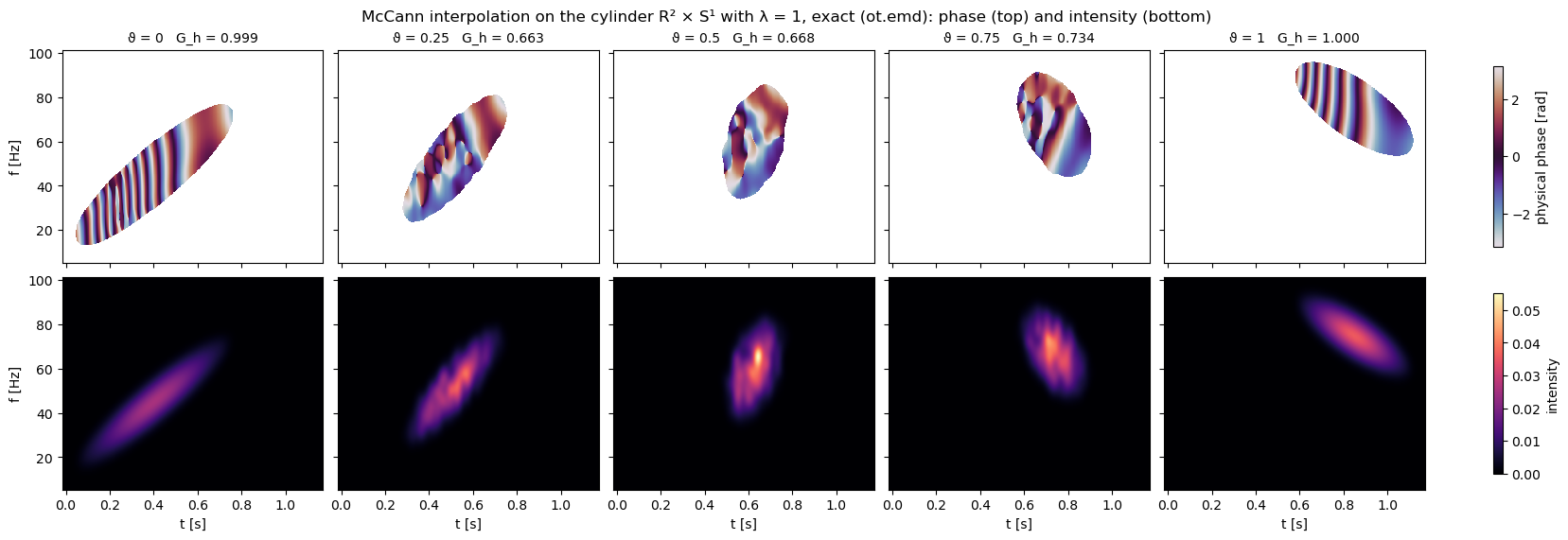}
    \caption{McCann interpolation with exact solver with respect to three geodesic distances: first row is the cc metric $(\cdot, \cdot)_{\Hb,\lambda}$; second row is the rescaled cylindrical Euclidean metric $g_\lambda\otimes d\theta^2$; and the third row is the cylindrical Euclidean metric $g_1\otimes d\theta^2$.}
    \label{fig:E1exact}
\end{figure}
\subsection{McCann interpolation along subRiemannian geodesics}\label{subsec: McCann}

Two linear chirps $f_0,f_1$, sampled at $f_s=512 \mathrm{Hz}$ over $1.25 \mathrm{s}$, are analyzed with a Gaussian window of standard deviation $\lambda=30 \mathrm{ms}$.
The STFT is evaluated on the uniform grid $t_n=n\Delta t$, $\omega_k=k\Delta\omega$ chosen so that $\Delta t/\lambda\approx\lambda\Delta\omega\approx0.5$. The hop $H$ and FFT length $N_{FFT}$, which satisfy $\Delta t=H/f_s$ and $\Delta\omega=2\pi f_s/N_{\mathrm{FFT}}$, are $8$ and $192$, respectively.
Using the soft cutoff $\varepsilon_{\mathrm{cut}}=0.05$, each retained coefficient defines an atom $(t_n,\omega_k,\Theta_f(t_n,\omega_k))\in\rH$, with normalized weight proportional to its cutoff amplitude.
The resulting endpoint lifts $\mu_0=\sum_i a_i\delta_{x_i},\mu_1=\sum_{j}b_j \delta_{y_j}$ consisting of $128$ atoms each are coupled using the quadratic cost $\mathring d_{\bH,\lambda}^2$, either exactly or with entropic regularization ($\varepsilon_{\mathrm{OT}}=0.25$).
For each transported pair, the interpolation follows the corresponding Carnot–Carathéodory geodesics $\gamma_{x_i\to y_j}$, so intermediate atoms generally leave the original STFT grid and no projection back onto that grid is performed.
Here $\gamma_{x_i\to y_j}$ is the projection to $\rH$ of $F_\lambda^{-1}\circ \gamma^\bG_{F_\lambda(x_i)\to F_\lambda(y_j)}$, the closed-form geodesic of Proposition \ref{prop:geodesics} pulled back by the isometry $F_\lambda$, where $y_j\star(0,0,2\pi k^*)$ is the nearest to $x_i$, obtained by wrapping the central coordinate of $F_\lambda(x_i^{-1}\star y_j)$ into $[-\pi,\pi]$ as in Proposition \ref{prop: wrapped geodesics}. 

\begin{figure}[!ht]
    \centering
    \includegraphics[width=1\linewidth]{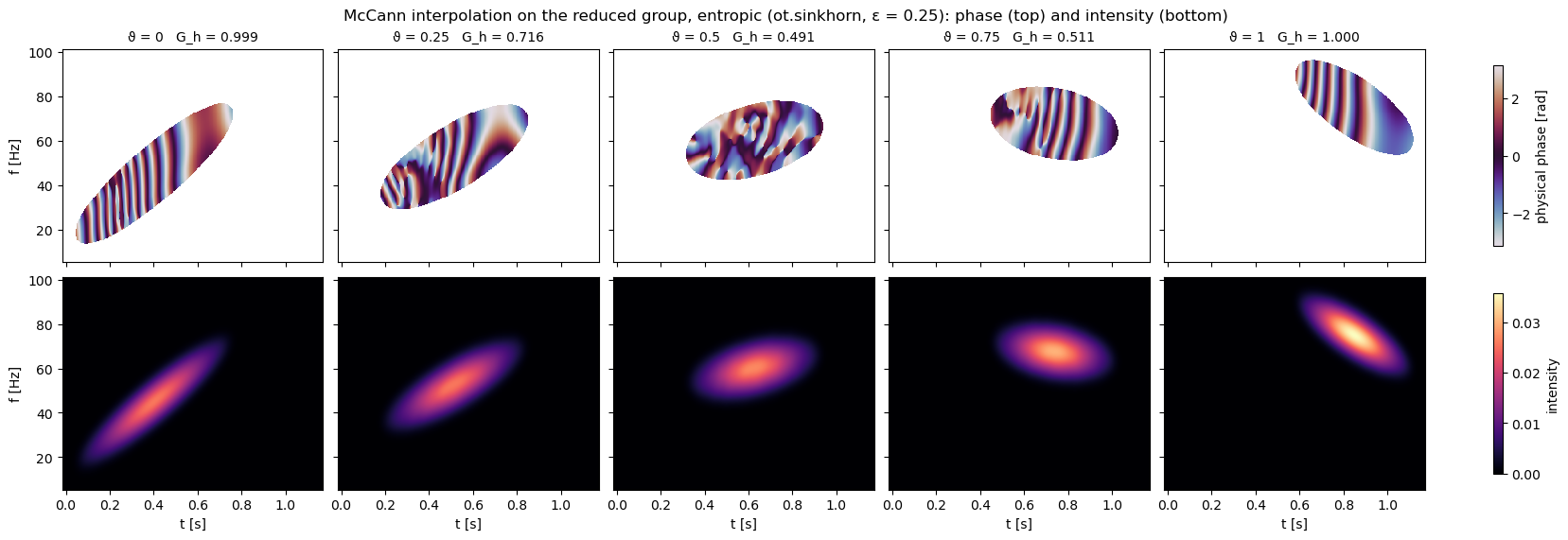}
    \includegraphics[width=1\linewidth]{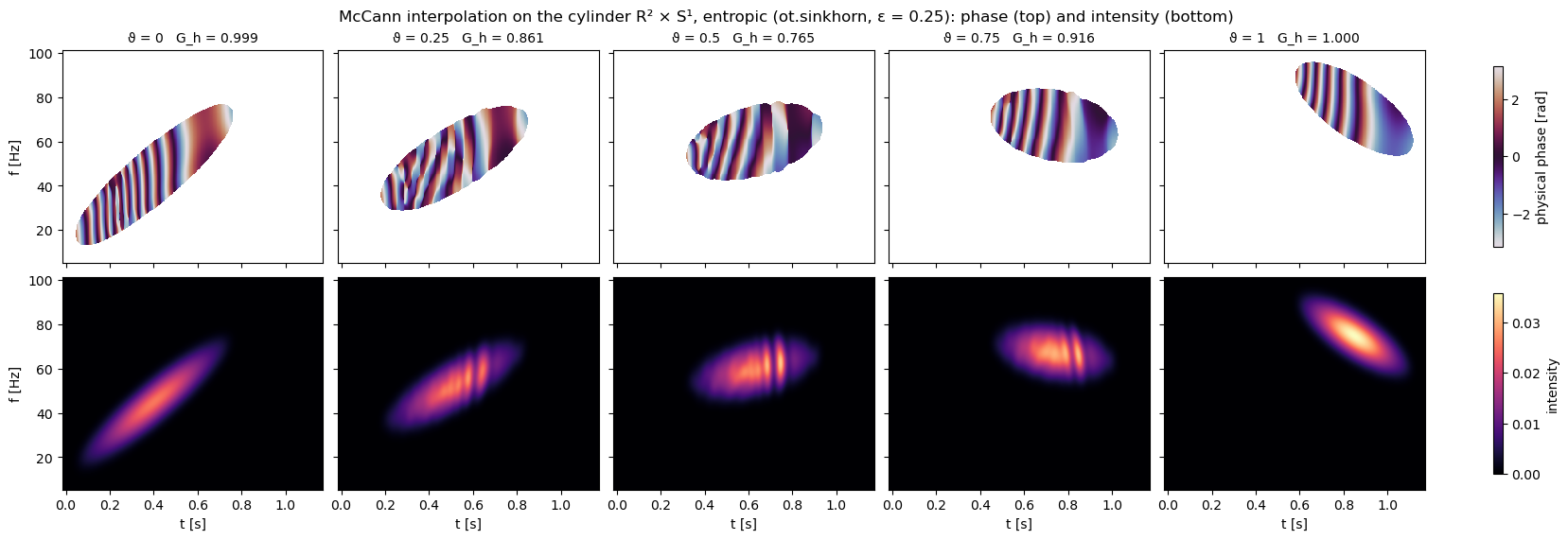}
    \caption{McCann interpolation with entropic regularization $\epsilon_{OT}=0.25$ with respect to two distances: first row is the Carnot-Carathéodory metric $(\cdot, \cdot)_{\Hb,\lambda}$; second row is the rescaled cylindrical Euclidean metric $g_\lambda\otimes d\theta^2$.}
    \label{fig:E1entr}
\end{figure}

In Figure \ref{fig:E1exact} and \ref{fig:E1entr}, an up-chirp around $45$ Hz at $0.40$ s and a down-chirp around $75$ Hz at $0.85$ s have disjoint lifts; 
The baseline is the flat cylinder $\mathbb R^2\times\mathbb T$ with the lifted metric $g_\lambda$ (Lemma \ref{lemma: projection}), whose geodesics are straight in $(t,\omega)$ and follow the shortest arc in $\theta$, for the same $\lambda$ and for $\lambda=1$ ($t$ in s, $\omega$ in rad/s). It shares the horizontal metric of $\mathring d_{\Hb,\lambda}$ but ignores the cocycle: its paths are not horizontal, and it is not left-invariant.

The Heisenberg interpolants move the energy in a single compact blob along an almost straight time–frequency path (first row of Figures \ref{fig:E1exact} and \ref{fig:E1entr}). 
In the rescaled Euclidean setting $g_\lambda$ (second row of Figures \ref{fig:E1exact} and \ref{fig:E1entr}), the interior intensity breaks into striations even at the entropy value $\epsilon_{OT}=0.25$; 
with $\lambda=1$ (third row of Figure \ref{fig:E1exact}) the frequency displacement dominates the cost and the chirp is compressed in time and fragmented. 
The Heisenberg phase is less regular to the eye: each atom carries the phase accumulated along its own horizontal geodesic, $\dot\theta=\omega\dot t$, and neighbours whose nearest lifts lie on opposite sides of $\sigma=\pm\pi$ can differ by about $\pi$ at $\vartheta=\frac12$. 
The Euclidean phase follows the shortest arc regardless of the motion in $(t,\omega)$. 

The graph index visible in Figures \ref{fig:E1exact} and \ref{fig:E1entr} confirms that the endpoint lifts are essentially phase graphs, with $\calG\simeq1$ in all cases. For the exact transport, the Heisenberg interpolants retain relatively high coherence, with $\calG=0.941,0.855,0.805$ at $\vartheta=0.25,0.5,0.75$; the scaled cylindrical metric gives comparable values $0.897,0.833,0.942$, whereas the unscaled choice $\lambda=1$ drops more markedly to $0.663,0.668,0.734$. Entropic regularization reduces phase coherence further: for the Heisenberg metric $\calG$ decreases to $0.716,0.491,0.511$, compared with $0.861,0.765,0.916$ for the scaled cylinder. Thus the graph index makes explicit the phase dispersion introduced by interpolation and, in particular, by entropic regularization.

\subsection{Convolutional Barycenters}\label{subsec: bary}
The interpolants of Section~\ref{subsec: McCann} were obtained by transporting atoms along closed-form geodesics, so that they leave the grid on which the lifts are computed. This subsection considers the complementary Eulerian setting, in which the lifts $\mu_{f_1},\dots,\mu_{f_n}$ of several signals are represented as weights on a single common grid of $\mathring{\mathbb H}$, and their entropic barycenters are computed on that same grid by the iterative Bregman projections of Section~\ref{subsec: computational}. The grid is built from the time-frequency sampling in the STFT in such a way that it becomes a finite window of the lattice $L_\Delta$ with $\Delta\phi = 2\pi/N_\phi$ and $N_\phi\,\Delta t\,\Delta\omega/2\pi \in \mathbb Z$, and is therefore contained in a discrete subgroup of $\mathring{\mathbb H}$. This compatibility with the group law is what makes the computation both exact and efficient. Every offset $z_k^{-1}\star z_l$ between grid points is again a lattice point, so the Gibbs kernel is sampled once on lattice offsets and each Sinkhorn update is an exact discrete group convolution, without any interpolation of the kernel or of the iterates. Moreover, since the cocycle of the group law is an integer multiple of $\Delta\phi$, a discrete Fourier transform in the periodic phase variable splits the group convolution into $N_\phi$ independent convolutions on the time-frequency grid: the zeroth mode is an ordinary planar convolution of the phase marginal, while every mode $m\neq 0$ is a twisted convolution, in which each time offset acts through an integer shear of the frequency-phase plane. The experiments below exploit this structure to compute barycenters of several lifted signals, to synthesize signals from them, and to verify that these barycenters follow a common time-frequency shift of the inputs exactly, a property not shared by the cylindrical baseline.

\paragraph{Lifting.}
Signals are lifted as in Definition~\ref{def: meas}, with a real Gaussian window of standard deviation $30$\,ms fixing $\lambda$ via \eqref{eq:lambda-choice}, the soft cutoff with $\varepsilon=0.05$, and the phase of each Gabor coefficient split linearly between its two neighbouring bins among $N_\phi=24$. The grid parameters are chosen so that the integer $N_\phi\Delta t\Delta\omega/2\pi$ is even. The half cocycle $t\omega/2$ is then also a whole number of phase bins, so that the symmetric coordinates are obtained by an exact relabelling of bins, and every lattice element $h$ acts on the discrete lift by an exact index shift, $\mu_{U_hf}=(L_h)_\#\mu_f$, away from the boundary of the grid.

\paragraph{Synthesis.}
The operator $R_\psi$ of Proposition~\ref{prop: rec} weights each time-frequency cell by the mean phasor $m(z)=\int_{\mathbb T}e^{i\theta}\,d\mu_z(\theta)$ of its fibre. This is exact on lifts, where $|m|=1$, but a barycenter spreads its mass along the fibres, so that only the fraction $\mathcal G(\mu)$ of the mass survives synthesis, and with phases dictated by the regularization. A barycenter does, however, reliably locate its mass in time and frequency, and by Proposition~\ref{prop:horizontal} this location determines the phase of a genuine signal through $d\phi=\omega\,dt$. Writing $b$ for a barycenter on the grid, its envelope and mean instantaneous frequency
\begin{equation}
a(t)=\sum_{\omega,\phi}b(t,\omega,\phi),\qquad
\bar\omega(t)=\frac{1}{a(t)}\sum_{\omega,\phi}\omega\,b(t,\omega,\phi),
\end{equation}
are therefore combined with the phase obtained by integrating the horizontality relation,
\begin{equation}
\phi^{\mathrm{syn}}(\tau)=\phi^\ast+\int_{t^\ast}^{\tau}\bar\omega(s)\,ds,
\qquad
f^{\mathrm{syn}}(\tau)=a(\tau)\,e^{i\phi^{\mathrm{syn}}(\tau)},
\end{equation}
where $(t^\ast,\omega^\ast,\phi^\ast)$ is the cell of largest mass and $a,\bar\omega$ are interpolated linearly from the frame times to the sampling times. The measured phase thus enters only through the global constant $\phi^\ast$, at the price of assuming a single component per time. A left translation by $h=(s,\eta,\theta)$ shifts $a$ by $s$, $\bar\omega$ by $\eta$ and the anchor to $h\star(t^\ast,\omega^\ast,\phi^\ast)$, so that $f^{\mathrm{syn}}\mapsto U_hf^{\mathrm{syn}}$ and the equivariance of the barycenters carries over to the synthesized signals.

\subsubsection{Multi-Marginal Entropic Barycenters}

\begin{figure}[t!]
  \centering
  \includegraphics[width=\textwidth]{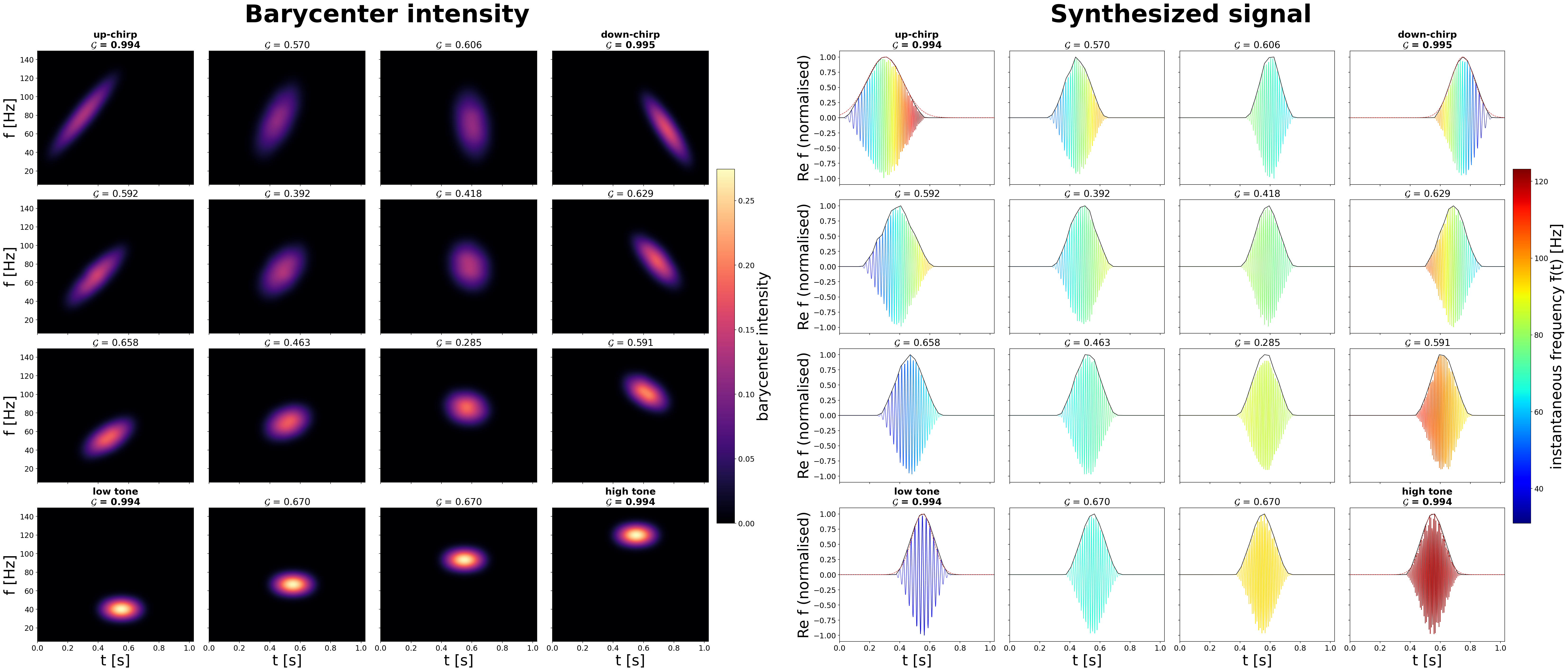}
  \caption{Entropic barycenters of four lifted signals (corners) with bilinear weights $\alpha$, Carnot-Carath\'eodory kernel, $\varepsilon_{\mathrm{OT}}=0.25$. \emph{Left:} time-frequency marginals. \emph{Right:} synthesized signals $\mathcal S(p)$, coloured by $\bar\omega(t)$, with the envelope $a(t)$ in black and, at the corners, the envelope of the input in dashed red. Each panel reports $\mathcal G$.}
  \label{fig:barycenter_grid}
\end{figure}

Figure~\ref{fig:barycenter_grid} shows entropic barycenters of four lifted signals, an up-chirp, a down-chirp, a low tone and a high tone, placed at the corners of a $4\times4$ array. The interior panels use the bilinear weights
\begin{equation}
\alpha(u,v)=\big((1-u)(1-v),\ u(1-v),\ (1-u)v,\ uv\big),\qquad u,v\in\left\{0,\tfrac13,\tfrac23,1\right\},
\end{equation}
so that every row and every column is a one-parameter family between two endpoints. Along the top row, the orientation of the up-chirp rotates through progressively shorter and more broadband atoms into that of the down-chirp, while along the bottom row the low tone is translated in frequency towards the high tone. In all cases the synthesized signals are clean single-component chirps and tones, whose slope and centre frequency vary continuously with $\alpha$, and the interior barycenters are shorter in time than the chirps, as expected from averaging long chirps with short tones. The graph index $\mathcal G$ equals $0.994$ to $0.995$ on the inputs, the ceiling $\cos(\Delta\phi/2)\approx0.991$ for a Dirac split between two phase bins being essentially attained, and drops to between $0.29$ and $0.67$ in the interior, with its lowest values where all four weights are nonzero. The entropic barycenters therefore spread substantially along the phase fibres, even for a pure frequency translation of a tone. This loss of phase coherence is invisible in the synthesized signals, since $\mathcal S$ reads only the time-frequency location of the mass and recovers the phase from horizontality.

\subsubsection{Equivariance to Time-Frequency Shifts}

\begin{figure}[b!]
  \centering
  \includegraphics[width=\textwidth]{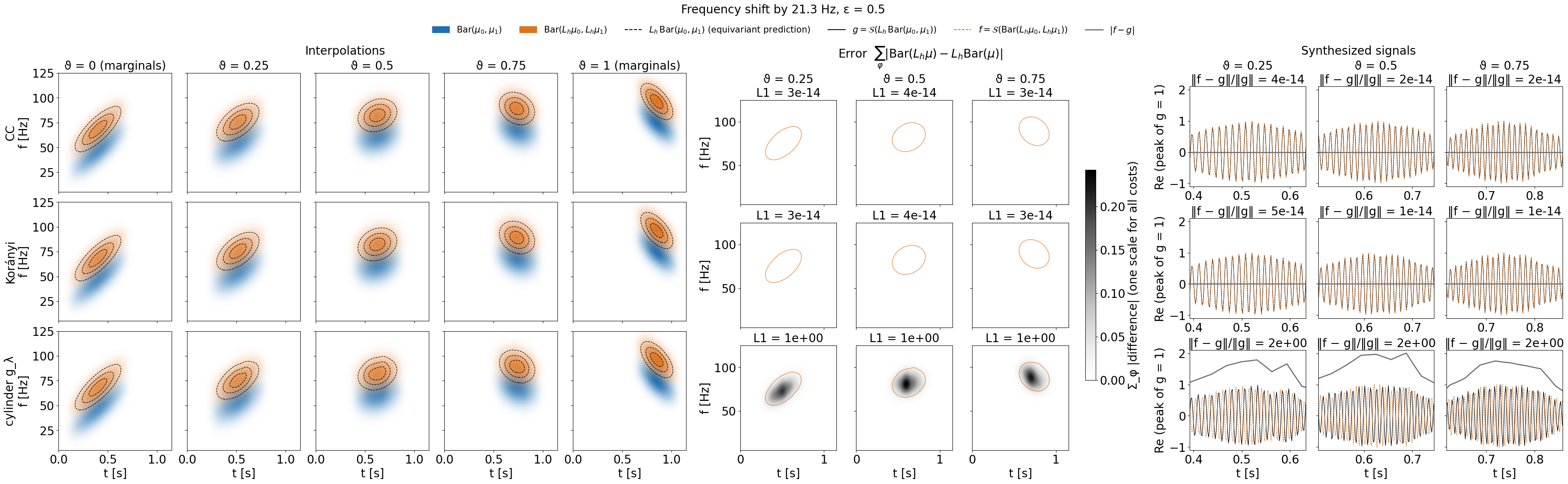}
  \caption{Equivariance under a frequency shift $h$ of $21.3$\,Hz, $\varepsilon_{\mathrm{OT}}=0.5$. Rows: Carnot-Carath\'eodory kernel, Kor\'anyi kernel, cylinder with metric $g_\lambda$. \emph{Left:} $\mathrm{Bar}(\mu_0,\mu_1)$ (blue), $\mathrm{Bar}(L_h\mu_0,L_h\mu_1)$ (orange) and $L_h\mathrm{Bar}(\mu_0,\mu_1)$ (dashed), shown through their time-frequency marginals. \emph{Middle:} $\sum_\phi|\mathrm{Bar}(L_h\mu)-L_h\mathrm{Bar}(\mu)|$ and its $L^1$ norm, on a common colour scale. \emph{Right:} $g=\mathcal S(L_h\mathrm{Bar}(\mu_0,\mu_1))$, $f=\mathcal S(\mathrm{Bar}(L_h\mu_0,L_h\mu_1))$ and $|f-g|$, with $\|f-g\|/\|g\|$.}
  \label{fig:equivariance_freq}
\end{figure}

Figure~\ref{fig:equivariance_freq} tests Remark~\ref{rmk: entropic-equivariance} on the grid. Two chirps are lifted to $\mu_0,\mu_1$, and the shifted inputs $L_h\mu_0,L_h\mu_1$ are obtained by lifting the shifted signals $U_hf_0,U_hf_1$, so that the comparison involves the whole pipeline. The barycenter $\mathrm{Bar}(L_h\mu_0,L_h\mu_1)$ with weights $(1-\vartheta,\vartheta)$ is then compared with the equivariant prediction $L_h\mathrm{Bar}(\mu_0,\mu_1)$. The shift $h=(0,\eta,0)$ is one frequency bin, $\eta/2\pi\approx21.3$\,Hz, and hence a lattice element. A frequency shift is the relevant test, since it acts on the phase by $\phi\mapsto\phi+\eta t$ and thus involves the cocycle of the group law, whereas a pure time shift leaves the phase coordinate unchanged. For both Heisenberg kernels the discrepancy is at round-off level, about $10^{-14}$ in $L^1$, for the measures as well as for the synthesized signals. This is expected, since the Gibbs kernels depend only on $z_k^{-1}\star z_l$, the counting measure on the lattice is invariant under lattice translations, and the discrete lift commutes with whole-bin shifts, so that every Bregman projection commutes with $L_h$. For the cylinder, the time-frequency marginals still nearly agree, because the planar part of its metric is translation invariant, but the phase-summed discrepancy is of order one. The map $\phi\mapsto\phi+\eta t$ is not an isometry of $\mathbb R^2\times\mathbb T$, since it changes the phase difference between atoms at different times by $\eta\,\Delta t$, so that the transport plan and the phase distribution of the barycenter change. In the synthesized signals, which share their envelope and instantaneous frequency, this appears as a global phase offset close to $\pi$, hence $\|f-g\|/\|g\|\approx2$.

\appendix
\section{The symmetrical Heisenberg group}\label{app: sym}
Let $\bG = \R^3_{(t,\omega,\sigma)}$ and define the differential 1-form, and codimension 1-tangent distribution
\begin{equation}
    \beta\coloneq d\sigma + \tfrac{1}{2}(t\,d\omega - \omega\,dt), \qquad E \coloneq \ker\beta.
\end{equation}

\begin{proposition}
\label{prop:symmetric}
The form $\beta$ is contact with Reeb field $R_\beta = \partial_\sigma$. A global horizontal frame is
\begin{equation}
  Y_1 = \partial_t + \frac{1}{2}\omega\partial_\sigma, \qquad
  Y_2 = \partial_\omega - \frac{1}{2}t\partial_\sigma,
\end{equation}
and $[Y_1, Y_2] = -\partial_\sigma = -R_{\beta}$. Declaring $Y_1, Y_2$ orthonormal gives a subRiemannian metric $(\cdot, \cdot)_\bG$ on $\bG$, and $Y_1,Y_2,R_{\beta},E,\beta$ are left-invariant for the group law
\begin{equation}\label{eq:heis-group}
  (t,\omega,\sigma) \odot (t',\omega',\sigma')
  \;=\;
  \Bigl(t+t',\; \omega+\omega',\; \sigma+\sigma'+\tfrac{1}{2}(\omega t'-t\omega')\Bigr).
\end{equation}
\end{proposition}

\begin{proof}
One has $d\beta = dt \wedge d\omega$, so $\beta \wedge d\beta \neq 0$.
The remaining claims follow by direct substitution and by differentiating left translations
for \eqref{eq:heis-group}.
\end{proof}
Note that the inverse of any $g \in \bG$ with respect to this group structure is simply $-g$. 
Analogously to the previous Section, we denote by $d_\bG$ the subRiemannian distance in $\bG$ associated with the metric $(\cdot, \cdot)_\bG$. In both cases, this distance is sometimes called Carnot-Carath\'eodory. With an abuse of notation, we still denote the actions $g\odot x$, and $x\odot g$ of a point $g\in \bG$ on $x\in \bG$ by $L_g, R_g\colon \bG\to \bG$, respectively.

\begin{proposition}\label{prop:geodesics}
    Every geodesic of $d_\bG$ from $0$ to some $x=(t,\omega, \sigma)\in \bG$ with $\sigma\neq 0$ is part of a curve of the type 
\begin{equation}\label{eq: geodesic}
    \gamma_{c,\theta}(\xi)=\left(e^{i \theta}\left(\frac{1- e^{-i\nu \xi}}{c}\right), \frac{\nu\xi -\sin(\nu \xi)}{2c^2}\right),
\end{equation}
for $c\in \R$ a curvature parameter, $\theta\in [0,2\pi]$ and $\nu=\operatorname{sgn}(c)\in \{\pm 1\}$ is the orientation. If $r^2=t^2+\omega^2>0$, the zero $s_0\in (-2\pi, 2\pi)$ of the function
\begin{equation}
    f(s)= \frac{r^2}{4} \frac{s-\sin (s)}{1-\cos(s)}- \sigma
\end{equation}
determines uniquely the geodesic by
\begin{equation}
    c=\operatorname{sgn}(s_0)\sqrt{\frac{2(1-\cos(s_0))}{r^2}}, \qquad \xi=|s_0|, \qquad \theta=\left[\frac{s_0}{2}-\arg(\omega+it)\right]\mod2\pi.
\end{equation}
In turn, $d_\bG(0, x)=s_0/c$.
\end{proposition}
\begin{proof}
    The first statement is precisely \cite[Equation 2.23]{Capogna}.
To identify the parameters $c,\theta$ and $\nu$ associated to the geodesic from the origin to $x$, we impose the passage from $x$ at time $\xi$ in $\gamma_{c,\theta}$:
\begin{equation}
    e^{i\theta}(1-e^{-i\nu \xi})=c(t+i\omega).
\end{equation}
Componentwise, we get
\begin{align}
    &ct=\cos(\theta)-\cos(\theta-\nu \xi)=2\sin\left(\frac{\nu \xi}{2}\right)\sin\left(\frac{\nu \xi}{2}-\theta\right);\\
    &c\omega=\sin(\theta)-\sin(\theta-\nu \xi)=2\sin\left(\frac{\nu \xi}{2}\right)\cos\left(\frac{\nu \xi}{2}-\theta\right);
\end{align}
By squaring and summing the two equations, we obtain the relation
\begin{equation}\label{eq: geo1}
    c^2(t^2+\omega^2)=4\sin^2\left(\frac{\nu \xi}{2}\right)=2(1-\cos(\nu \xi)).
\end{equation}
Regarding the third component, we assume $r^2\coloneq t^2+\omega^2> 0$ and use Equation \eqref{eq: geo1} to write
\begin{align}
    \sigma=\frac{\nu\xi -\sin(\nu \xi)}{2c^2}= \frac{r^2}{4} \frac{\nu\xi-\sin (\nu \xi)}{1-\cos(\nu \xi)}.
\end{align}
Substituting $s=\nu \xi$, we get the equation is satisfied at the zero $s_0\in (-2\pi,2\pi)$ of $f$.
\end{proof}
The formula of Proposition \ref{prop:geodesics} parametrizes the nonstraight geodesics with nonzero horizontal radius. The case $c=0$ is obtained as a limit and gives horizontal straight lines. 
Points on the center, $t=\omega=0$, lie in the cut locus of the origin and admit nonunique minimizing geodesics.

The Carnot-Carath\'eodory distance is intrinsic but requires solving a geodesic problem. For estimates and numerical kernels it is often useful to work with an algebraic homogeneous gauge.
\begin{definition}\label{def:koranyi}
    We define the Kor\'anyi gauge on $\bG$ as
    \begin{equation}
    |(t,\omega,\sigma)|_K\coloneq((t^2+\omega^2)^2+16\sigma^2)^{1/4}.
    \end{equation}
    Consequently, we define the Kor\'anyi distance on $\bG$ as
    \begin{align}
        & d_K(x, y)\coloneq |y^{-1}\odot x|_K=|L^{-1}_y x|_K.
    \end{align}
\end{definition}
The Kor\'anyi distance is naturally left-invariant and preserves the anisotropic scaling of the group. Moreover, it is bi-Lipschitz equivalent to the Carnot-Carath\'eodory distance (cf Proposition \ref{prop:cc-koranyi}), while avoiding the root-finding step required by the exact geodesic distance. 

\begin{proposition}\label{prop:cc-koranyi}
For all $x,y\in \bG$,
\begin{equation}\label{eq:cc-koranyi}
d_K(x,y)\;\le\; d_{\bG}(x,y)\;\le\; \sqrt{\pi}\,d_K(x,y),
\end{equation}
with equality on the left if $y^{-1}\odot x\in\{\sigma=0\}$ and on the right if $y^{-1}\odot x\in\{t=\omega=0\}$. 
\begin{proof}
    The bounds follow from the structure of CC geodesics, see \cite[Equation 1.4]{Balogh2006} and \cite{Bella97}.
\end{proof}
\end{proposition}
\begin{remark}
    The choice of the coefficient $1/2$ for the contact and group structure in $\bG$ is by convention. Substituting $1/2$ by any $\zeta\neq0$ would change the geodesic equation \eqref{eq: geodesic} in its central component ($\zeta<0$ would flip the orientation), making it scale as $2\zeta$ times the enclosed area.
 \end{remark}

\begin{proposition}\label{prop:equi}
Define the map
\begin{equation}
  F_\lambda\colon \Hb \longrightarrow \bG, \qquad
  F_\lambda(t,\omega,\phi) \;=\; \Bigl(\frac{t}{\lambda},\; \lambda \omega,\; \phi - \tfrac{1}{2}t\omega\Bigr).
\end{equation}
Then
\begin{equation}
  (F_\lambda)^*\beta = \alpha, \qquad
  (F_\lambda)_*(X_1^\lambda) = Y_1, \qquad (F_\lambda)_*(X_2^\lambda) = Y_2,
\end{equation}
and $F_\lambda(h \star h') = F_\lambda(h) \odot F_\lambda(h')$.
Therefore $F_\lambda$ is simultaneously a contactomorphism, a Carnot-group isomorphism,
and a global subRiemannian isometry.
\end{proposition}

\begin{proof}
Since $d\sigma = d\phi - \tfrac{1}{2}(\omega\,dt + t\,d\omega)$, we get
$d\sigma + \tfrac{1}{2}(t\,d\omega - \omega\,dt) = d\phi - \omega\,dt$.
The frame identities follow by differentiating $F_\lambda$.
Using \eqref{eq:jet-group},
\begin{equation}
  \phi + \phi' + \omega t' - \tfrac{1}{2}(t+t')(\omega+\omega')
  = \Bigl(\phi - \tfrac{1}{2}t\omega\Bigr) + \Bigl(\phi' - \tfrac{1}{2}t'\omega'\Bigr)
    + \tfrac{1}{2}(\omega t' - t\omega'),
\end{equation}
which is exactly the third component of \eqref{eq:heis-group} applied to $F_\lambda(h)$ and $F_\lambda(h')$.
\end{proof}
Proposition \ref{prop:equi} allows us to transfer every metric statement between the two coordinate systems. In particular,
\begin{equation}
d_{\Hb,\lambda}(h,h')=d_\bG(F_\lambda(h),F_\lambda(h')),
\end{equation}
and a minimizing geodesic from $h$ to $h'$ in $\Hb$ is obtained by applying $F_\lambda^{-1}$ to a minimizing geodesic from $F_\lambda(h)$ to $F_\lambda(h')$ in $G$. We also transport the Kor\'anyi gauge to the original coordinates by 
\begin{equation}
    |h|_{K,\lambda}\coloneq|F_\lambda(h)|_K,\qquad d_{K,\lambda}(h,h')\coloneq d_K(F_\lambda(h),F_\lambda(h')).
\end{equation}
In particular, \eqref{eq:cc-koranyi} holds verbatim for $d_{K,\lambda}$ and $d_{\Hb,\lambda}$ on $\Hb$.



\paragraph{Declaration of generative AI-assisted technologies in the manuscript preparation.}
The idea of adapting the computational optimal transport framework of \cite{pai2025} to the Heisenberg setting in harmonic analysis is due to the authors, after a few failed attempts to apply the Heisenberg group in different contexts. The lifted measure construction without the cutoff was also part of the initial idea, although the model \textrm{GPT-5.6} later suggested the introduction of the cutoff to deal with phase singularities near zero-amplitude points. The same model also suggested the use of a window-dependent parameter $\lambda$ in the Carnot–Carathéodory metric, which was later observed to be sensible in applications. The proofs contained in the manuscript are original and were written by the authors, although they were checked and refined by the aforementioned model or the following version \textrm{GPT-6 Sol}. Technologies such as Claude Code (\textrm{Opus 5.5}) and OpenCode (\textrm{GPT-6.1 Sol}) were used to assist with the numerics. Both models were also used to review the final edit of the manuscript in terms of narrative and correctness, and \textrm{GPT-6.1 Sol} suggested the application of our framework to delay-Doppler synchronization in radar signal processing, as described in Remark \ref{rmk: radar}. All outputs have been thoroughly checked by the authors and were obtained through several interactions with the models. The authors take full responsibility for the content.

\small
\bibliographystyle{plain}
\bibliography{heisenberg}

\end{document}